\documentclass[a4article,12pt,reqno]{amsart}
\usepackage{amsmath,amsfonts,amsthm,amssymb}
\usepackage{epsfig}

\usepackage{amsmath,amsfonts,amsthm,amssymb,color}
\usepackage[latin1]{inputenc}

\newcommand{\ist}{\int_{s}^{t}}
 
\newcommand{\ott}{[0,T]}

\newcommand{\id}{\mbox{id}}

\newcommand{\lot}{[\ell_1,\ell_2]}
\newcommand{\norm}[1]{\lVert #1\rVert}
\newcommand{\xd}{\mathbf{x}^{\mathbf{2}}}

\newcommand{\xti}{\tilde{x}}
\newcommand{\yti}{\tilde{y}}
\newcommand{\zti}{\tilde{z}}

\newcommand{\bb}{\mathbf{B}}
\newcommand{\bd}{{\bf B^{2}}}
\newcommand{\bt}{\mathbf{2}}
\newcommand{\xdt}{{\bf \tilde{x}^{2}}}
\newcommand{\der}{\delta}
\newcommand{\II}{{\rm i}}

\newcommand{\E}{\mathbf{E}}
\newcommand{\cZ}{{\mathcal Z}}

\def\half{{\frac{1}{2}}}

\definecolor{rpp}{rgb}{0.5, 0.5, 1.0}
\definecolor{rp}{rgb}{0.25, 0, 0.75}
\definecolor{dg}{rgb}{0, 0.5, 0}

\newcommand{\C}{\mathbb C}

\newcommand{\R}{\mathbb R}
\newcommand{\N}{\mathbb N}

\newcommand{\cb}{\mathcal B}
\newcommand{\cac}{\mathcal C}

\newcommand{\cd}{\mathcal D}

\newcommand{\cj}{\mathcal J}

\newcommand{\cn}{\mathcal N}
\newcommand{\cq}{\mathcal Q}

\newcommand{\al}{\alpha}

\newcommand{\ga}{\gamma}
\newcommand{\gga}{\Gamma}

\newcommand{\laa}{\Lambda}

\newcommand{\si}{\sigma}

\newcommand{\vp}{\varphi}

\newcommand{\ka}{\kappa}

\newcommand{\lp}{\left(}
\newcommand{\rp}{\right)}
\newcommand{\lc}{\left[}
\newcommand{\rc}{\right]}
\newcommand{\lcl}{\left\{}
\newcommand{\rcl}{\right\}}
\newcommand{\lln}{\left|}
\newcommand{\rrn}{\right|}

\newcommand{\bean}{\begin{eqnarray*}}
\newcommand{\eean}{\end{eqnarray*}}
\newcommand{\ben}{\begin{enumerate}}
\newcommand{\een}{\end{enumerate}}
\newcommand{\beq}{\begin{equation}}
\newcommand{\eeq}{\end{equation}}

\newtheorem{theorem}{Theorem}[section]

\newtheorem{corollary}[theorem]{Corollary}

\newtheorem{definition}[theorem]{Definition}

\newtheorem{hypothesis}{Hypothesis}
\newtheorem{lemma}[theorem]{Lemma}

\newtheorem{proposition}[theorem]{Proposition}

\theoremstyle{remark}
\newtheorem{remark}[theorem]{Remark}

\begin{document}
\title[A Milstein-type scheme  for SDEs driven by FBM]{A Milstein-type scheme without L\'evy area terms for SDEs driven by fractional Brownian motion
}
\author{A. Deya, A. Neuenkirch, S. Tindel}
\date{\today}
\begin{abstract}
In this article, we study the numerical approximation of stochastic differential equations driven by a multidimensional fractional Brownian motion (fBm) with Hurst parameter greater than $1/3$.  We introduce an implementable scheme for these equations, which is based on a second order Taylor expansion, where the usual L\'evy area terms are replaced  by products of increments of the driving  fBm. The convergence of our scheme is shown by means of a combination of rough paths techniques and error bounds for the discretisation of the L\'evy area terms.
\end{abstract}

\address{Aur{\'e}lien Deya, Institut {\'E}lie Cartan Nancy, Universit\'e de Nancy 1, B.P. 239,
54506 Vand{\oe}uvre-l{\`e}s-Nancy Cedex, France}
\email{deya@iecn.u-nancy.fr}

\address{Andreas  Neuenkirch, Technische Universit\"at Dortmund, Fakult\"at f\"ur Mathematik, 
Vogelpothsweg 87,
D-44227 Dortmund,
Germany} \email{andreas.neuenkirch@math.tu-dortmund.de}

\address{Samy Tindel, Institut {\'E}lie Cartan Nancy, Universit\'e de Nancy 1, B.P. 239,
54506 Vand{\oe}uvre-l{\`e}s-Nancy Cedex, France}
\email{tindel@iecn.u-nancy.fr}

\thanks{Supported by the DAAD (PPP-Procope D/0707564) and Egide (PHC-Procope 17879TH)}

\date{\today}

\maketitle

\section{Introduction and Main Results}

Fractional Brownian motion (fBm in short for the remainder of the article) is a natural generalisation of the usual Brownian motion, insofar as it is defined as a centered Gaussian process $B=\{B_t;\, t \in \R_+ \}$ with continuous sample paths, whose increments $(\der B)_{st}:=B_t-B_s$, $s,t \in \R_+$ are characterised by their variance $\E[(\der B)_{st}^2]=|t-s|^{2H}$. Here  the parameter $H \in (0,1)$, which is called Hurst parameter,  governs  in particular the H\"older regularity of the sample paths of $B$ by a standard application of Kolmogorov's criterion: fBm has H\"older continuous sample paths of order $\lambda$ for all $\lambda < H$. 
 The particular case $H=1/2$ corresponds to the usual Brownian motion,  so the cases $H\ne 1/2$ are a natural extension of the classical situation, allowing e.g. any prescribed H\"older regularity of the driving process. Moreover, fBm is $H$-self similar, i.e. for any $c>0$ the process
$\{ c^H B_{t/c}; \, t \in \R_+ \}$
is again a fBm, and also has stationarity increments, that is for any $h \geq 0$ the process
$ \{ B_{t+h}-B_h; \, t \in \R_+ \}$
is a fBm.

\smallskip

These properties  (partially) explain why stochastic equations driven by fBm have received considerable attention during the last two decades. Indeed, many physical systems seem to be governed by a Gaussian noise with different properties than  classical Brownian motion.
Fractional Brownian motion as driving noise is used e.g. in electrical engineering \cite{DMS, DR}, or biophysics \cite{CLCHH,Ko, OTHP}. Moreover, after some controversial discussions (see  \cite{bjork} for a summary of the early developments) fBm has established itself also in financial modelling, see e.g. \cite{Guasoni, Bender}. For empirical studies of fractional Brownian motion in finance see e.g.  \cite{CKW, W2, Corcuera}.
 All these situations lead to different kind of stochastic differential equations (SDEs), whose simplest prototype can be formally written as 
\beq\label{eq:frac-sde}
Y_t= a+ \sum_{i=1}^m\int_0^t \sigma^{(i)}(Y_u) \, d B_u^{(i)}, \quad t \in \lc 0, T \rc, \qquad  a \in \R^d,
\eeq
where  $\si=(\si^{(1)},\ldots,\si^{(m)})$ is a smooth enough function  from $\R^d$ to $\R^{d\times m}$ and $B=(B^{(1)},\ldots,B^{(m)})$ is a $m$-dimensional fBm with Hurst parameter $H>1/3$.

\smallskip

At a mathematical level, fractional differential equations of type (\ref{eq:frac-sde}) are typically handled (for $H\ne 1/2$) by pathwise or semi-pathwise methods. Indeed for $H>1/2$, the integrals $\int_0^t \sigma^{(i)}(Y_u) \, d B_u^{(i)}$, $i=1, \ldots, m,$  in (\ref{eq:frac-sde}) can be defined  using Young integration or fractional calculus tools, and these methods also yield the  existence of a unique solution,
see e.g. \cite{NR,Za}. When $1/4<H<1/2$, the existence and uniqueness result for equation~(\ref{eq:frac-sde}) can be seen as the canonical example of an application of the rough paths theory. The reader is referred to \cite{FV,Lyons-bk} for the original version of the rough paths theory, and to \cite{Gu} for a (slightly) simpler algebraic setting which will be used in the current article. In the particular case $1/3<H<1/2$, the rough path machinery can be summarised  very briefly as follows: assume that our driving signal $B$ allows to define iterated integrals with respect to itself. Then one can define and solve equation (\ref{eq:frac-sde}) in a reasonable class of processes.

\smallskip

Once SDEs driven by  fBm are solved, it is quite natural (as in the case of SDEs driven by the usual Brownian motion) to study the stochastic processes they define. However, even if some progress has been made in this direction, e.g.  concerning the law of the solution \cite{BC,CFV,NNRT} or its ergodic properties \cite{HO}, the picture here is far from being complete. Moreover, explicit solutions of stochastic differential equations driven by fBm are rarely known, as in the case of SDEs  driven by classical Brownian motion. Thus one has to rely on numerical methods for the simulation of these equations.

\smallskip

So far,  some numerical schemes for equations like (\ref{eq:frac-sde}) have already been studied in the literature. In the following,  we consider uniform grids of the form $\{t_k=kT/n;\, 0\le k \le n\}$ for a fixed $T>0$. The simplest  approximation method is the Euler scheme defined by 
\begin{align*}
Y^n_{0}&=a, \\
{Y}_{t_{k+1}}^{n} &= {Y}_{t_k}^{n}   +  
\sum_{i=1}^m \sigma^{(i)} ({Y}_{t_{k}}^{n}) \der B^{(i)}_{t_{k}t_{k+1}}, \qquad k=0, \ldots, n-1.
\end{align*}
For $H>1/2$, the Euler scheme converges to the solution of the SDE (\ref{eq:frac-sde}). See e.g. in~\cite{MS}, where an almost sure  convergence  rate $n^{-(2H-1)+ \varepsilon}$ with $\varepsilon >0$ arbitrarily small  is established.   A detailed analysis of the one-dimensional case is given in  \cite{NN}, where 
the exact convergence rate  $n^{-2H+1}$ and  the asymptotic error distribution are derived.

\smallskip

However, the Euler scheme is not appropriate  to approximate SDEs driven by fBm when $1/3<H<1/2$. This is easily illustrated by the following one-dimensional example, in which $B$ denotes a one-dimensional fBm: consider the equation $$dY_t=Y_t \, dB_t,\quad t \in [0,1], \qquad Y_0=1,$$
whose exact solution is  $$Y_t=\exp(B_t), \quad t \in [0,1].$$ The Euler approximation for this equation at the final time point $t =1$ can be written as
$$  Y^{n}_1 = \prod_{k=0}^{n-1} (1+ (\der B)_{k/n,(k+1)/n}).$$
So for $n \in \mathbb{N}$ sufficiently large and using a Taylor expansion, we have
$$
Y^{n}_1=\exp\Big( \sum_{k=0}^{n-1} \log  (1+ (\der B)_{k/n,(k+1)/n})\Big)
=\exp\Big (B_1 -\frac12 \sum_{k=0}^{n-1}|(\der B)_{k/n,(k+1)/n}|^2 + \rho_n \Big),
$$
where  $\lim_{ n \rightarrow \infty} \rho_n \stackrel{\rm a.s.}{=} 0 $ for $H>1/3$. Now it is well known that 
$$ \sum_{k=0}^{n-1}|(\der B)_{k/n,(k+1)/n}|^2 \stackrel{\rm a.s.}{ \longrightarrow} \infty $$ for $H<1/2$
as $n\to\infty$, which implies that  $  \lim_{ n \rightarrow \infty} Y^{n}_1 \stackrel{{\rm a.s.}}{=} 0 $. This is obviously incompatible with a convergence towards $Y_1=\exp(B_1)$.  In the case $H=1/2$ this phenomenon is also well known: here the Euler scheme converges to the It\^o solution and not to the Stratonovich solution of SDE (\ref{eq:frac-sde}).

\smallskip

To obtain a convergent numerical method Davie proposed in \cite{Da} a scheme of Milstein type. For this, assume that all  iterated integrals of $ B$ with respect to itself are collected into a $m\times m$ matrix  $\bb^\bt$, i.e. set
$$ \bb^\bt_{st}(i,j) = \int_s^t  \int_s^u dB_{v}^{(i)} \,  dB_{u}^{(j)}, \qquad 0\le s<t\le T, \quad 1 \le i,j\le m.$$ 
The matrix $\bb^\bt$ (respectively its elements) is (are) usually called L\'evy area.
Davie's scheme is then given by
\begin{align}\label{eq:davie-scheme} 
Y^n_{0}&=a, \\
{Y}_{t_{k+1}}^{n} &= {Y}_{t_k}^{n}   +  
\sum_{i=1}^m \sigma^{(i)} ({Y}_{t_{k}}^{n})\der B^{(i)}_{t_{k}t_{k+1}}+
 \sum_{i,j=1}^m \mathcal{D}^{(i)}\sigma^{(j)} ( Y_{t_{k}}^{n}) \, \bb^\bt_{t_{k}t_{k+1}}(i,j),\nonumber  \qquad k=0, \ldots, n-1,\end{align}
with the  differential operator $\mathcal{D}^{(i)} = \sum_{l=1}^{d} \sigma_l^{(i)} \partial_{x_l}$. (Recall that we use the notation $\der B^{(i)}_{st}=B^{(i)}_{t}-B^{(i)}_{s}$ for $s,t \in \lc 0, T \rc$.)  This scheme is shown to be convergent as long as $H>1/3$ in \cite{Da}, with an almost sure convergence rate of $n^{-(3H-1)+\varepsilon}$ for $\varepsilon >0$ arbitrarily small. This result has then been extended in \cite{FV} to an abstract rough path with arbitrary regularity, under further assumptions on the higher order iterated integral of the driving signal. 

\smallskip

As the classical Milstein scheme for SDEs driven by Brownian motion, the Milstein-type scheme  (\ref{eq:davie-scheme}) is in general not a directly implementable method. Indeed, unless
the commutativity condition $$   \mathcal{D}^{(i)}\sigma^{(j)} =   \mathcal{D}^{(j)}\sigma^{(i)}, \qquad i,j=1, \ldots, m,$$ holds,
the simulation of the iterated  integrals $\bb^\bt_{t_{k}t_{k+1}}(i,j)$ is necessary. However, the law of these integrals is unknown, so that they can not be simulated directly and have to be approximated.

 In this article we replace the iterated integrals by a simple product of increments, i.e. we use the approximation
\begin{align} \label{disc_levy} \bb^\bt_{t_{k}t_{k+1}}(i,j)  \approx \frac{1}{2}  \, \der B^{(i)}_{t_{k}t_{k+1}} \, \der B^{(j)}_{t_{k}t_{k+1}}. \end{align}
This leads to the following simpler Milstein-type scheme:
Set $Z_{t_{0}}^{n}=a$ and 
 \begin{align} \label{eq:misltein-scheme}
Z_{t_{k+1}}^{n} =  Z_{t_k}^{n}   
+  \sum_{i=1}^m \sigma^{(i)} ( Z_{t_{k}}^{n}) \, \der B^{(i)}_{t_{k}t_{k+1}}
+  \frac{1}{2} \sum_{i,j=1}^m \mathcal{D}^{(i)}\sigma^{(j)} ( Z_{t_{k}}^{n})
\, \der B^{(i)}_{t_{k}t_{k+1}} \, \der B^{(j)}_{t_{k}t_{k+1}}
\end{align}
for $ k=0, \ldots, n-1$.
Moreover, for $t \in (t_k,t_{k+1})$, define
 \begin{align} \label{eq:misltein-scheme_2} Z_{t}^{n} = Z_{t_k}^{n} + \frac{t-t_k}{T/n} \big (  \der Z^n \big)_{t_kt_{k+1}}, \end{align}
i.e. if $t \in [0,T]$ is not a discretisation point, then   $Z_{t}^{n}$ is defined by piecewise linear interpolation. This scheme  is now directly implementable and is still convergent.

\bigskip

\begin{theorem} \label{main_thm} Assume that $\si \in C^3(\R^d; \R^{d\times m})$ is bounded with bounded derivatives. Let $Y$ be the  solution to equation (\ref{eq:frac-sde}) and $Z^{n}$  the Milstein approximation given by (\ref{eq:misltein-scheme}) and~(\ref{eq:misltein-scheme_2}). Moreover, let $1/3 < \gamma  < H$. Then, there exists a finite and non-negative random variable $\eta_{H ,\gamma, \si , T}$ such that
\begin{equation}\label{main_thm_eq}
    \| Y- Z^{n} \|_{\gamma, \infty,T}   \leq  \eta_{H, \gamma, \si , T} \cdot  \sqrt{\log(n)} \cdot n^{-(H-\gamma)}
\end{equation} for $n>1$.
\end{theorem}

\bigskip

Here $\| \cdot \|_{\kappa, \infty,T}$ denotes the $\kappa$-H\"older norm of a function $f: 
\lc 0, T \rc \rightarrow \R^l$, i.e.
\begin{align}
\|f \|_{\kappa, \infty,T} = \sup_{t \in \lc 0 ,T \rc} |f(t)| + \sup_{s,t \in \lc 0 ,T \rc} \frac{|f(t)-f(s)|}{|t-s|^{\kappa}}. \label{hoelder-norm}
\end{align}

\begin{remark}
Note that the almost sure estimate (\ref{main_thm_eq}) cannot be turned into an $L^1$-estimate for $\| Y-Z^n\|_{\ga,\infty,T}$. This is a common consequence of the use of the rough paths method, which exhibits non-integrable (random) constants, as a careful examination of the proof of Theorem \ref{thm:Lip} would show. See also \cite{FV} for further details.
\end{remark}

\smallskip

Our strategy to prove the above Theorem consists of two steps.
First we determine the error between $Y$ and its Wong-Zakai approximation
\beq\label{WZ}
\overline{Z}^n_t= a+ \sum_{i=1}^m\int_0^t \sigma^{(i)}(\overline{Z}_u^n) \, d B_u^{(i),n}, \quad t \in \lc 0, T \rc, \qquad  a \in \R^d,
\eeq
where
$$
B^{n}_t=B_{t_k}+\lp \frac{t-t_k}{T/n}\rp (\der B)_{t_{k}t_{k+1}}, \quad  t \in \lc 0, T \rc,
$$ 
i.e. $B$ in equation (\ref{eq:frac-sde}) is replaced with its piecewise linear interpolation. (For a survey on Wong-Zakai approximations for standard SDEs see e.g. \cite{Wong-Zakai}.)
Here, we denote the  L\'evy area corresponding to $B^n$ by  $\bb^n$.
 Using the Lipschitzness of the It\^o map of $Y$, i.e. the solution of equation (\ref{eq:frac-sde}) depends continuously in appropriate H\"older norms   on $B$ and the L\'evy-area $\bb$,  and error bounds for the difference between $B$ and $B^n$ resp. $\bb$ and $\bb^n$, we obtain
  $$ \| Y- \overline{Z}^{n} \|_{\gamma, \infty,T}   \leq   \eta_{H ,\gamma, \si , T}^{(1)} \cdot  \sqrt{\log(n)} \cdot n^{-(H-\gamma)},
$$
 where $ \eta_{H ,\gamma, \si , T}^{(1)}$ is  a finite and non-negative random variable.\\
 In the second step we analyse the difference between $\overline{Z}^n$ and $Z^n$.
 The second order Taylor scheme with stepsize $T/n$ for classical ordinary differential equations applied to  the Wong-Zakai approximation
(\ref{WZ})  gives our simplified Milstein scheme (\ref{eq:misltein-scheme}). So to obtain the  error bound   $$      \| Z^n- \overline{Z}^{n} \|_{\gamma, \infty,T}   \leq   \eta_{H ,\gamma, \si , T}^{(2)} \cdot  \sqrt{\log(n)} \cdot n^{-(H-\gamma)},
$$
we can proceed in a similar way as for the numerical analysis of classical ordinary differential equations.
We first determine the one-step error  and then control the error propagation using  a global stability result with respect to the initial value for differential equations driven by rough paths. The latter can be considered as a substitute for Gronwall's lemma in this context.

\medskip

Combining both error bounds then gives Theorem \ref{main_thm}.

\medskip

\begin{remark}\label{class-Milstein}  For $H=1/2$ the scheme 
(\ref{eq:davie-scheme})  corresponds to the classical Milstein scheme   for Stratonovich SDEs driven by Brownian motion, while our scheme (\ref{eq:misltein-scheme}) corresponds to the so called simplified Milstein scheme. See e.g. \cite{KP}. 
\end{remark}

\begin{remark}
At the price of further computations, which are simpler than the ones in this article,  our convergence result can be extended to an equation with drift, i.e. to
$$
Y_t= a
+ \int_0^t b(Y_u) \, du
+ \sum_{i=1}^m\int_0^t \sigma^{(i)}(Y_u) \, d B_u^{(i)}, \quad t \in \lc 0, T \rc, \qquad  a \in \R^d,
$$
where $b:\R^d\to\R^d$ is a $C_b^3$ function and where the other coefficients satisfy the assumptions of Theorem \ref{main_thm}. Indeed, the equation above can be treated like our original system~(\ref{eq:frac-sde}) by adding a component $B^{(0)}_t=t$ to the fractional Brownian motion. The additional iterated integrals of $B^{(0)}$ with respect to $B^{(j)}$ for $j=1,\ldots,m$ are easier to handle  than $\bb^\bt(i,j)$ for $i,j \in \{1, \ldots, m\}$, since they are classical Riemann-Stieltjes integrals.  For sake of conciseness we do not include the  corresponding details.
\end{remark}

\begin{remark}\label{non_bounded} Theorem \ref{main_thm} requires $ \si$ to be bounded. However, if $\si \in C^3(\R^d; \R^{d\times m})$ is neither bounded nor has bounded derivatives but equation (\ref{eq:frac-sde})  has still a unique pathwise solution in the sense of Theorem \ref{thm:Lip} below, then the assertion of Theorem \ref{main_thm} is still valid. This follows from a standard localisation procedure, see e.g.  \cite{JKN},  and   applies in particular to affine-linear coefficients.
\end{remark}

\begin{remark}\label{rem_low_bound} The error bound of Theorem \ref{main_thm} is sharp. To see this, consider the most simple equation 
$$dY^{(1)}_t=dB^{(1)}_t, \quad t\in \lc 0,T \rc, \qquad Y_0=a\in \R,$$
for which our  approximation obviously reduces to $Z^n=B^n$. Then, due to results of H\"usler, Piterbarg and Seleznjev (\cite{huesleretal}) for the deviation of a Gaussian process from its linear approximation, one can prove that
\begin{equation}\label{eq:norm-sup-hold}
\lim_{n \rightarrow \infty}      {\bf P} \left(\,  \ell(n) \cdot \| Y- Z^{n} \|_{\gamma, \infty,T} \, < \, \infty \right)=0  , $$
 if $$\liminf_{n \rightarrow \infty}  \, \, {\ell(n)} \cdot  \sqrt{\log(n)} \cdot n^{-(H-\gamma)} = \infty .
\end{equation}
For further details see  Section \ref{optimality}. 
\end{remark}

\begin{remark}\label{other_disc}   If the Wong-Zakai approximation is discretised with an arbitrary numerical scheme  for ODEs of at least second order (e.g. Heun, Runge-Kutta 4), then the arising scheme for equation (\ref{eq:frac-sde}) satisfies the same error bound as the proposed modified Milstein scheme. So, the strategy of our proof is in fact an instruction for the construction of arbitrary implementable and  convergent numerical schemes for SDEs driven by fBm.
\end{remark}

\begin{remark}\label{alt_approach}  Instead of  replacing  the  L\'evy terms in Davie's scheme
by the "rough" approximation  (\ref{disc_levy}) one could discretise these terms very finely using the results contained 
 in \cite{NTU}, where (exact) convergence rates for approximations of the L\'evy area are derived. However, it is well known that already for SDEs driven by Brownian motion such a scheme is rarely efficient, if the convergence rate of the scheme is measured in terms of its computational cost. For a survey on the complexity  of the approximation of SDEs driven by Brownian motion, see e.g.
 \cite{thomas}.
\end{remark}

\medskip

The $\gamma$-H\"older norm, which appears  in Theorem \ref{main_thm}  since the It\^o-map of $Y$ is only Lipschitz in appropriate H\"older norms with $1/3 < \gamma <H$  and thus is natural in the rough path setting, is not typical for measuring the error of  approximations to stochastic differential equations.  A more standard criterion would be  the  error with respect to the supremum norm, i.e.
$$ \| Y - Z^{n} \|_{\infty, T} = \sup_{t \in \lc 0, T \rc }  | Y_t - Z^{n}_t |.$$
The error (in the supremum norm) of the piecewise linear interpolation of fractional Brownian motion is of order $\sqrt{\log(n)} n^{-H}$, see \cite{huesleretal}. Moreover, for the iterated integral
$ \int_0^T \int_0^u dB_{v}^{(1)} \,  dB_{u}^{(2)}$ the proposed Milstein-type scheme leads to the trapezoidal type approximation
$$ \frac{1}{2}\sum_{k=0}^{n-1}  \big( B^{(1)}_{t_k} + B^{(1)}_{t_{k+1}}\big) \big( B^{(2)}_{t_{k+1}} - B^{(2)}_{t_{k}} \big)  .$$
The $L^p$-error  for this approximation is of order $n^{-2H+1/2}$, see \cite{NTU}.

Based on these two findings, our guess for the rate of convergence in supremum norm is that
$$ \| Y - Z^{n} \|_{\infty, T}  \leq \eta_{H, \sigma, T} \cdot  \sqrt{\log(n)} \cdot \big( n^{-H} + n^{-2H+1/2} \big)$$
{holds under the assumptions of Theorem \ref{main_thm}. 
This conjecture is also supported by the numerical examples we give in Section 4.

\smallskip

The remainder of this article is structured as follows: In Section \ref{sec:alg-integration} we recall some basic facts on algebraic integration and rough differential equations. The proofs of Theorem~\ref{main_thm} and Remark \ref{rem_low_bound}  are given in Section 3 and 4.  Finally, Section 5  contains the mentioned numerical examples.

\section{Algebraic integration and differential equations}
\label{sec:alg-integration}

In this section, we recall the main concepts of algebraic integration, which will be essential  to define  the generalized integrals in our setting. Namely, we state the definition of the spaces of increments, of the operator $\delta$, and its inverse called $\Lambda$ (or sewing map according to the terminology of \cite{FP}). We also recall some elementary but useful algebraic relations on the spaces of increments. The interested reader is sent to \cite{Gu} for a complete account on the topic, or to \cite{DT,GT} for a more detailed summary.

\subsection{Increments}\label{sec:incr}

The extended integral we deal
with is based on the notion of increments, together with an
elementary operator $\der$ acting on them. 

The notion of increment can be introduced in the following way:  for two arbitrary real numbers
$\ell_2>\ell_1\ge 0$, a vector space $V$, and an integer $k\ge 1$, we denote by
$\cac_k( [\ell_1,\ell_2]; V)$ the set of continuous functions $g : [\ell_1,\ell_2]^{k} \to V$ such
that $g_{t_1 \cdots t_{k}} = 0$
whenever $t_i = t_{i+1}$ for some $i\in \{0, \ldots,  k-1\}$.
Such a function will be called a
\emph{$(k-1)$-increment}, and we will
set $\cac_*( [\ell_1,\ell_2];V)=\cup_{k\ge 1}\cac_k( [\ell_1,\ell_2]; V)$. 
To simplify the notation, we will write $\cac_k(V)$, if there is no ambiguity about $ [\ell_1,\ell_2] $.

The operator $\der$
is  an operator acting on
$k$-increments, 
and is defined as follows on $\cac_k(V)$:
\begin{equation}
  \label{eq:coboundary}
\delta : \cac_k(V) \to \cac_{k+1}(V), \qquad
(\delta g)_{t_1 \cdots t_{k+1}} = \sum_{i=1}^{k+1} (-1)^i g_{t_1
  \cdots \hat t_i \cdots t_{k+1}} ,
\end{equation}
where $\hat t_i$ means that this particular argument is omitted.
Then a fundamental property of $\der$, which is easily verified,
is that
$\delta \delta = 0$, where $\delta \delta$ is considered as an operator
from $\cac_k(V)$ to $\cac_{k+2}(V)$.
 We will denote $\cZ\cac_k(V) = \cac_k(V) \cap \text{Ker}\delta$
and $\cb \cac_k(V) =\cac_k(V) \cap \text{Im}\delta$.

\smallskip

Some simple examples of actions of $\der$,
which will be the ones we will really use throughout the article,
 are obtained by letting
$g\in\cac_1(V)$ and $h\in\cac_2(V)$. Then, for any $t,u,s\in\lot$, we have
\begin{equation}
\label{eq:simple_application}
  (\der g)_{st} = g_t - g_s
\quad\mbox{ and }\quad
(\der h)_{sut} = h_{st}-h_{su}-h_{ut}.
\end{equation}

\smallskip

Our future discussions will mainly rely on
$k$-increments with $k =2$ or $k=3$, for which we will use some
analytical assumptions. Namely,
we measure the size of these increments by H\"older norms
defined in the following way: for $f \in \cac_2(V)$ let
$$
\norm{f}_{\mu} = \sup_{s,t\in\lot}\frac{|f_{st}|}{|t-s|^\mu}
\quad\mbox{and}\quad
\cac_2^\mu(V)=\lcl f \in \cac_2(V);\, \norm{f}_{\mu}<\infty  \rcl.
$$
Using this notation, we define in a natural way 
$$
\cac_1^\mu(V)=\{ f \in \cac_1(V); \, \norm{\der f}_\mu < \infty \},
$$
and recall that we have also defined a norm $\|\cdot\|_{\kappa,\infty,T}$ at equation (\ref{hoelder-norm}).
In the same way, for $h \in \cac_3(V)$, we set
\begin{eqnarray}
  \label{eq:normOCC2}
  \norm{h}_{\gamma,\rho} &=& \sup_{s,u,t\in\lot}
\frac{|h_{sut}|}{|u-s|^\gamma |t-u|^\rho},\\
\norm{h}_\mu &= &
\inf\left \{\sum_i \norm{h_i}_{\rho_i,\mu-\rho_i} ;\, h  =\sum_i h_i,\, 0 < \rho_i < \mu \right\} ,\nonumber
\end{eqnarray}
where the last infimum is taken over all sequences $\{h_i, \, i \in \mathbb{N}\} \subset \cac_3(V) $ such that $h
= \sum_i h_i$ and over all choices of the numbers $\rho_i \in (0,\mu)$.
Then  $\norm{\cdot}_\mu$ is easily seen to be a norm on $\cac_3(V)$, and we define
$$
\cac_3^\mu(V):=\lcl h\in\cac_3(V);\, \norm{h}_\mu<\infty \rcl.
$$
Eventually,
let $\cac_3^{1+}(V) = \cup_{\mu > 1} \cac_3^\mu(V)$,
and note  that the same kind of norms can be considered on the
spaces $\cZ \cac_3(V)$, leading to the definition of the  spaces
$\cZ \cac_3^\mu(V)$ and $\cZ \cac_3^{1+}(V)$. In order to avoid ambiguities, we denote  in the following  by $\cn[  \, \cdot \,  ; \cac_j^\kappa]$ the $\kappa$-H\"older norm on the space $\cac_j$, for $j=1,2,3$. For $\zeta\in\cac_j(V)$, we also set $\mathcal{N}[\zeta;\mathcal{C}_{j}^{0}(V)]=\sup_{s\in[\ell_1; \ell_2]^j}\lVert \zeta_s\rVert_{V}$.

\smallskip

The operator $\der$ can be inverted under some H\"older regularity conditions, which is essential for the construction of our generalized integrals. 

\begin{theorem}[The sewing map] \label{prop:Lambda}
Let $\mu >1$. For any $h\in \cZ \cac_3^\mu( V)$, there exists a unique $\Lambda h \in \cac_2^\mu(V)$ such that $\der( \Lambda h )=h$. Furthermore,
\begin{eqnarray} \label{contraction}
\norm{ \Lambda h}_\mu \leq \frac{1}{2-2^{\mu}}\, \cn[h;\, \cac_3^{\mu}(V)].
\end{eqnarray}
 This gives rise to a  continuous linear map $\laa:  \cZ \cac_3^\mu( V) \rightarrow \cac_2^\mu(V)$ such that $\der \laa =\id_{ \cZ \cac_3^\mu( V)}$.
\end{theorem}

\begin{proof}
The original proof of this result can be found in \cite{Gu}. We refer to \cite{DT,GT} for two simplified versions.

\end{proof}

The sewing map creates a first link between the structures we just introduced and the problem of integration of irregular functions:

\begin{corollary}[Integration of small increments]
\label{cor:integration}
For any 1-increment $g\in\cac_2 (V)$ such that $\der g\in\cac_3^{1+}$,
set
$
h = (\id-\Lambda \delta) g
$.
Then, there exists $f \in \cac_1(V)$ such that $h=\der f$ and
$$
(\delta f)_{st} = \lim_{|\Pi_{st}| \to 0} \sum_{i=0}^n g_{t_{i} t_{i+1}},
$$
where the limit is over any partition $\Pi_{st} = \{t_0=s,\dots,
t_n=t\}$ of $[s,t]$ whose mesh tends to zero. The
1-increment $\delta f$ is the indefinite integral of the 1-increment $g$.
\end{corollary}

\smallskip

We also need some product rules for the operator $\delta$. For this
recall the following convention:
for  $g\in\cac_n( \lot \R^{l,d})$ and $h\in\cac_m( \lot ;\R^{d,p}) $ let  $gh$
be the element of $\cac_{n+m-1}( \lot ;\R^{l,p})$ defined by
\begin{equation}\label{cvpdt}
(gh)_{t_1,\dots,t_{m+n-1}}=g_{t_1,\dots,t_{n}} h_{t_{n},\dots,t_{m+n-1}}
\end{equation}
for $t_1,\dots,t_{m+n-1}\in \lot$. With this notation, the following elementary rule holds true:

\begin{proposition}\label{difrul}
Let $g\in\cac_2 (\lot;\R^{l,d})$ and $h\in\cac_1(\lot;\R^d)$. Then
$gh$ is an element of $\cac_2( \lot ;\R^l)$ and
$\der (gh) = \der g\, h  - g \,\der h.$
\end{proposition}

\smallskip

\subsection{Random differential equations}
One of the main appeals of the algebraic integration theory is that 
differential equations driven by   a $\ga$-H\"older signal $x$ can be defined and solved rather quickly in this setting. In the case of an H\"older exponent $\ga>1/3$, the required structures  are just the notion of \textit{controlled processes} and the L\'evy area based on $x$.

\smallskip

Indeed, let us consider an equation of the form
\begin{align} \label{de}
dy_t=\sigma(y_t) \, dx_t=\sum_{i=1}^{m}\sigma^{(i)}(y_t) \, dx_t^{i}, 
\quad t \in \lc 0, T \rc, \qquad  y_0=a,
\end{align}
where $a$ is a given initial condition in $\R^d$, $x$ is an element of $\cac_1^{\ga}([0,T];\, \R^m)$, and $\si$ is a smooth enough function from $\R^d$ to $\R^{d,m}$. Then it is  natural (see \cite{TU} for further explanations) that the increments of a candidate  for a solution to (\ref{de}) should be controlled by the increments of $x$ in the following way:
\begin{definition}
\label{def:controlled-process}
Let $z$ be a path in $\cac_1^\ka(\R^d)$ with $1/3<\ka\le\ga$.
We say that $z$ is a weakly controlled path based on $x$ if
$z_0=a$, with $a\in\R^d$,
and $\der z\in\cac_2^\ka(\R^d)$ has a decomposition $\der z=\zeta \der x+ r$,
that is, for any $s, t\in[0,T]$,
\begin{equation}
\label{f3.0.4}
(\der z)_{st}=\zeta_s (\der x)_{st} + r_{st},
\end{equation}
with $\zeta\in\cac_1^\ka(\R^{d,m})$ and $r\in\cac_2^{2\ka}(\R^d)$. 
\end{definition}

\noindent
The space of weakly controlled
paths will be denoted by $\cq^x_{\ka,a}(\R^d)$, and a process
$z\in\cq^x_{\ka,a}(\R^d)$ can be considered in fact as a couple
$(z,\zeta)$.  The space $\cq^x_{\ka,a}(\R^d)$ is endowed with a natural semi-norm given by
\begin{align}
\label{f3.0.5}
& \cn[z;\cq^x_{\ka,a}(\R^d)]  \\
 & \qquad = \cn[z;\cac_1^{\ka}(\R^d)]
+ \cn[\zeta;\cac_1^{0}(\R^{d,m})]
+ \cn[\zeta;\cac_1^{\ka}(\R^{d,m})]
+\cn[r;\cac_2^{2\ka}(\R^d)], \nonumber
\end{align}
where the quantities $\cn[g;\cac_j^{\ka}]$ have been defined in Section \ref{sec:incr}.
For the L\'evy area associated to $x$ we assume the following structure:
\begin{hypothesis}\label{h1}
The path $x:[0,T]\to\R^m$ is 
$\ga$-H\"older continuous with $\frac{1}{3}<\ga\le 1$ and  admits a so-called L\'evy area,
that is, a process 
$\xd\in\cac_2^{2\ga}(\R^{m,m})$, which  satisfies
$\der\xd=\der x\otimes \der x$, namely
\begin{equation*}
\lc (\der\xd)_{sut} \rc(i,j)=[\der x^{i}]_{su} [\der x^{j}]_{ut},
\end{equation*}
for any $s,u,t\in[0,T]$ and $i,j\in\{1,\ldots,m  \}$. 
\end{hypothesis}

\smallskip

To illustrate the idea behind the construction of the generalized  integral 
assume that the paths $x$ and $z$ are smooth  and also for simplicity that $d=m=1$. Then the Riemann-Stieltjes integral of $z$ with respect to $x$ is well defined and 
we have
$$
\ist z_udx_u = z_s(x_t-x_s) + \ist(z_u-z_s)dx_u  = z_s( \der x)_{st} + \ist (\der z)_{su}dx_u 
$$
for $\ell_1 \leq s\le t \leq \ell_2 $.   If $z$ admits the decomposition (\ref{f3.0.4})
we
obtain
\begin{align}
\ist (\der z)_{su}dx_u  =
\ist \left( \zeta_{s} (\der x)_{su}  +\rho_{su} \right)  dx_u
=  \zeta_{s} \ist (\der x)_{su}  \, dx_u + \ist \rho_{su}  \, dx_u.
\label{heuristic}
\end{align}
Moreover, if we set
$$( \xd)_{st}:= \ist (\der x)_{su}  \, dx_u, \qquad \ell_1 \leq s\le t \leq \ell_2,$$
then it is quickly verified that $\xd$ is the associated L\'evy area to $x$. Hence we can write
$$
\ist  z_udx_u 
=  z_s (\der x)_{sz} +  \zeta_{s}\,  (\xd)_{st}  + \ist \rho_{su}  \, dx_u.
$$ 
 Now rewrite this equation as
\begin{align}
  \ist \rho_{su}  \, dx_u = \ist  z_udx_u  - z_s (\der x)_{st} -  \zeta_{s}\,  (\xd)_{st}
\label{exp1:irhodx}
\end{align}
and apply  the increment operator  $\der$ to both  sides of this equation.
For smooth paths $z$ and $x$
we have
$$
\der \left (\int z \, dx \right)=0,
\qquad \qquad
\der(z \,\der x)= - \der z \, \der x,
$$
by Proposition \ref{difrul}.
Hence, applying these relations to the right hand side of
(\ref{exp1:irhodx}), using the decomposition (\ref{f3.0.4}), the properties of the L\'evy area and
again Proposition \ref{difrul},
we obtain
\begin{align*}
 \left[\der \left( \int \rho \, dx \right) \right]_{sut}  &= \, 
(\der z)_{su} (\der x)_{ut}
+ (\der\zeta)_{su} \, (\xd)_{ut}
-\zeta_{s}  \, (\delta \xd )_{sut} \\
&= \zeta_s   (\der x)_{su}   \, (\der x)_{ut}
+\rho_{su} \, (\der x)_{ut}
+ (\der\zeta)_{su} (\xd)_{ut}
-\zeta_{s}  (\der x)_{su}   \, (\der x)_{ut}
\\
&=\rho_{su} (\der x)_{ut}+  (\der\zeta )_{su} \, (\xd)_{ut}.
\end{align*}
So in summary, we have derived  the representation
$$
\left[\der \left( \int \rho \, dx \right) \right]_{sut}  =\rho_{su} (\der x)_{ut}+  (\der\zeta )_{su} \, (\xd)_{ut}.
$$
As we are dealing with smooth paths we have
$
\der \left( \int \rho \, dx \right) \in\mathcal{Z}\cac_3^{1+}$ and thus  belongs to the domain
 of $\laa$   due  to
Proposition \ref{prop:Lambda}.  (Recall that $\der\der=0$.) Hence, it follows
$$
\ist \rho_{su} \, dx_u =
\laa_{st}\lp \rho\, \der x  + 
\der\zeta \, \xd \rp,
$$ and
 inserting this identity into (\ref{heuristic}) we end up with
$$
\ist  z_udx_u 
=  z_s (\der x)_{st} +  \zeta_{s}\,  (\xd)_{st}  + 
\laa_{st}\lp \rho\, \der x  + 
\der\zeta \, \xd \rp.
$$
Since in addition
$$     \rho\, \der x  + 
\der \zeta \, \xd =   - \der(z \der x + \zeta \, \xd)  ,$$
we can also write this as
$$ \int  z_udx_u = (\id-\Lambda \delta) (z \der x + \zeta \, \xd).$$
Thus we have expressed the Riemann-Stieltjes integral of $z$ with respect to $x$
in  terms of the sewing map $\laa$, of the  L\'evy area $\mathbf{x}^{\mathbf{2}}$
 and of increments of $z$ resp. $x$. This can now be generalized to the non-smooth case.
Note that Corollary \ref{cor:integration} justifies the use of the notion integral.

In the following, we denote by $A^*$ the transposition of a vector resp. matrix,
and by
 $A_1 \cdot A_2 ={\rm Tr}(A_1A_2^*)$ the inner product of two vectors or
two matrices $A_1$ and $A_2$.

\begin{proposition}\label{intg:mdx}
For fixed $\frac13 < \ka \leq \ga$,
let $x$ be a path satisfying Hypothesis \ref{h1}. Furthermore,  let
$z \in\cq^x_{\ka,\alpha}([\ell_1,\ell_2];\R^{m})$ such that the
increments of $z$ are given by (\ref{f3.0.4}).
Define $\hat{z}$ by $\hat{z}_{\ell_1}= \hat{\alpha}$ with $\hat{\alpha} \in \R$  and
\begin{equation}\label{dcp:mdx}
(\der \hat{z})_{st}= \lc (\operatorname{id}-\Lambda \delta) (z^* \der x + \zeta \cdot   \xd) \rc_{st}
\end{equation}
for $\ell_1 \leq s \leq t \leq \ell_2$.
Then $\cj(z^*\, dx):= \hat{z} $ is a well-defined element of $\cq^x_{\ka, \hat{\alpha}}([\ell_1,\ell_2];\R)$ and coincides with the usual Riemann integral, whenever
$z$ and $x$ are smooth functions. 
\end{proposition}

Moreover, the H\"older norm of $\cj(z^*\, dx)$ can be estimated in terms of the H\"older norm of the integrator $z$. (For this and also for a proof of the above Proposition, see e.g.  \cite{Gu}.) This allows to use a fixed point argument to obtain the existence of a unique solution for rough differential equations.

\begin{theorem}\label{thm:Lip}
For fixed $\frac13 < \ka < \ga$, let $x$ be a path satisfying Hypothesis \ref{h1}, and let
$\si \in C^3(\R^{d}; \R^{d,m})$ be bounded with bounded derivatives. Then we have:
\begin{enumerate}
\item
Equation (\ref{de}) admits a unique solution $y$  in
$\cq^x_{\ka,a}([0,T];\R^d)$ for  any $T>0$, and there exists a polynomial $P_T: \R^2 \rightarrow \R^+$ such that
\begin{equation}\label{control-sol-y}
\cn[y;\cq^x_{\ka,a}([0,T];\R^d)] \leq P_T(\| x\|_{\gamma,\infty,T},\|  \xd \|_{2 \gamma})
\end{equation}
holds.
\item
 Let
$F: \R^d \times \cac_1^{\ga}([0,T];\R^m)\times
\cac_2^{2\ga}([0,T];\R^{m, m}) \rightarrow
  \cac_1^{\ga}([0,T];\R^d)$  be the mapping defined by $$ F\left( a ,x,\xd  \right) =y, $$
where $y$ is the unique solution of equation (\ref{de}).  This mapping
is locally
 Lipschitz continuous in the following sense:
Let $\tilde{x} $ be another driving rough path with corresponding
 L\'evy area $\xdt $ and $\tilde{a}$ be another initial condition. Moreover denote by
$ \tilde{y}$ the unique solution of the corresponding differential equation.
        Then, there exists an increasing function $K_T: \R^4 \rightarrow \R^+$
such that
\begin{align}\label{lipsch-cont}
  \qquad \,\,\,  \| y - \tilde{y}\|_{\ga,\infty,T}  
  &\leq K_T(\| x\|_{\gamma,\infty,T} ,  \| \tilde{x} \|_{\gamma,\infty,T},  \|  \xd \|_{2 \gamma}, \|  \xdt \|_{2 \gamma}  ) \, \\& \qquad \qquad \times    \left(  |a - \tilde{a}| +  \| x - \tilde{x}\|_{\gamma,\infty,T} +  
\|  \xd- \xdt \|_{2 \gamma}  \right) \nonumber
\end{align}
holds,
where we recall that $\|f\|_{\mu,\infty,T}=\|f\|_{\infty}+\| \der f \|_{\mu}$ denotes  the usual H\"older norm of a path
$f \in \cac_{1}([0,T]; \R^l)$.
\end{enumerate}
\end{theorem}

\begin{remark}
{\it Inequality (\ref{control-sol-y}) implies in particular
\begin{equation}\label{estim-remainder-term}
\lln (\der y)_{st}-\si(y_s) (\der y)_{st} \rrn \leq \lln t-s \rrn^{2\ka}  P_T(\| x\|_{\gamma,\infty,T},\|  \xd \|_{2 \gamma}).
\end{equation}
This estimate will be required in the proof of Lemma \ref{cont-wz}.}
\end{remark}

The above Theorem improves (slightly) the original formulation of the Lipschitz continuity of the It\^o map $F$, which can be found in \cite{Gu}, concerning the control of the solution in terms of the driving signal. Therefore (and also for completeness) we provide some details of its proof  in the appendix.  A similar continuity result can be found in \cite{FV}, where the classical approach  of  Lyons and Qian  to rough differential equations is used.

\bigskip

\subsection{Application to fBm}
The application of the rough path theory to an equation with a particular driving signal relies on the existence of the L\'evy area fulfilling Hypothesis~\ref{h1}. In our setting, the driving process is given by an $m$-dimensional fractional Brownian motion
$(B^{(1)}, \ldots, B^{(m)})$ with Hurst parameter $\gamma >1/3$.

To the best of our knowledge,  there are three known possibilities to show the existence of the associated L\'evy area $\bb^{\bt}=(\bb^{\bt}(i,j))_{i,j=1, \ldots, m}$:
 (i) By a piecewise dyadic linear interpolation of the paths of $B$, as done in  \cite{CQ}. (ii) Using Malliavin calculus tools in order to define $\bb^\bt$ as a Russo-Vallois iterated integral, similarly to what is done in~\cite{NNT} to construct a delayed fractional L\'evy area. (iii) By means of the analytic approximation of $B$ introduced by Unterberger in~\cite{Un}.
Actually, all three methods lead to the same L\'evy area. The equivalence between the first two constructions has been established by Coutin and Qian through a representation formula (see Theorem 4 in \cite{CQ}). The convergence results we are going to establish show that the L\'evy area recently obtained by Unterberger in \cite{Un} coincide with the previous ones. Note that this question had been left open by the author in the latter reference, so that the following Proposition \ref{main_wz_Lp} has an interest in itself (see also \cite{NTU} for a partial result in this direction).

\smallskip

We resort here to the analytic definition of the fractional L\'evy area, since we  use the pointwise estimates of \cite{NTU}, which were derived in this setting. Let us recall the main features of the analytic approach.

\smallskip

\subsubsection{Definition of the analytic fBm}
The article \cite{Un} introduces the fractional Brownian motion as the real part of 
the trace on $\R$ of an analytic process $\gga$ (called: {\em analytic fractional
Brownian motion} \cite{TU}) defined on the complex upper-half plane $\Pi^+=\{z\in\C;\, \Im(z)>0\}$. This is achieved by an   explicit series construction: for $k\ge 0$ and $z\in\Pi^+$, set
\begin{equation} 
f_k(z)=2^{H-1} \sqrt{\frac{H(1-2H)}{2\cos\pi H}}
 \sqrt\frac{\mathbf{\Gamma}(2-2H+k)}{\mathbf{\Gamma}(2-2H) k!}
\left( \frac{z+\II}{2\II} \right)^{2H-2} \left( \frac{z-\II}{z+\II}\right)^k,  
\end{equation}
where $\mathbf{\Gamma}$ stands for the usual Gamma function. These functions are well-defined on $\Pi^+$, and it can be checked that 
\begin{equation*} 
\sum_{k\ge 0} f_k\lp x+\II \frac{\eta_1}{2}\rp \overline{f_k\lp y+\II\frac{\eta_2}{2}\rp} = K^{',-}\lp\half\lp\eta_1+\eta_2\rp;x,y\rp,
\end{equation*} 
where $K^{',-}$ is a positive kernel defined on $\R_+^*\times\R\times\R$ given by
$$
K^{',-}(\eta;x,y)=\frac{H(1-2H)}{2\cos\pi H} (-\II(x-y)+\eta)^{2H-2}.
$$
We also set
$$
K^{',+}(\eta;x,y)=\frac{H(1-2H)}{2\cos\pi H} (+\II(x-y)+\eta)^{2H-2}.
$$

\smallskip

Now define the  Gaussian process $\Gamma'$ with ''time parameter'' $z\in \Pi^+$ by
\begin{equation} \label{eq:7}
\Gamma'(z)=\sum_{k\ge 0} f_k(z)\xi_k 
\end{equation}
 where $(\xi_k)_{k\ge 0}$ are independent standard complex
Gaussian variables, i.e. $\E[\xi_j \xi_k]=0$, $\E[\xi_j \bar{\xi}_k]=\delta_{j,k}$.
The Cayley transform $z\mapsto  \frac{z-\II}{z+\II}$ maps $\Pi^+$ to $\cd$, where $\cd$ stands for the unit disk of the complex plane. This allows to prove  that the series defining $\Gamma'$ is a random entire series which is analytic on the unit disk and hence the process $\Gamma'$ is analytic on $\Pi^+$. Furthermore, restricting to the horizontal line $\R+\II\frac{\eta}{2}$, the following identity holds:
$$ 
\E[\Gamma'(x+\II\eta/2)\overline{\Gamma'(y+\II\eta/2)}]=K^{',-}(\eta;x,y).
$$
 
\smallskip

One may now integrate the process $\Gamma'$ over any path $\gamma:(0,1)\to\Pi^+$ with endpoints $\gamma(0)=0$
and $\gamma(1)=z\in\Pi^+\cup\R$ (the result does not depend on the particular path but only on the endpoint $z$).
The resulting  process, which is denoted by  $\Gamma$,  is still analytic on $\Pi^+$. Furthermore,
 the real part of the boundary value of $\Gamma$
on $\R$ is a fractional Brownian motion. Another way to look at this is to define $\Gamma(\eta):= \{ \Gamma(t+\II\eta); t \in \R\}$
 as a regular process living on $\R$, and to observe that the real part of $\Gamma(\eta)$ converges
for $\eta\to 0$ to a fractional Brownian motion.  The following Proposition summarises what has been said so far:

\begin{proposition}[see \cite{Un,TU}]
Let $\gga'$ be the process defined on $\Pi^+$ by relation (\ref{eq:7}).
\begin{enumerate}
\item
Let $\gamma:(0,1)\to\Pi^+$ be a continuous path with endpoints $\gamma(0)=0$ and $\gamma(1)=z$,
 and set $\Gamma_z=\int_{\gamma}\Gamma'_u \, du$. Then $\Gamma$ is an analytic process on $\Pi^+$.
 Furthermore, as $z$ runs along any path in $\Pi^+$ going to $t\in\R$, the random variables $\Gamma_z$
 converge almost surely to a random variable called again $\Gamma_t$.

\item

The family $\{\Gamma_t;\, t\in\R\}$ defines a  centered Gaussian complex-valued process whose
paths are almost surely $\kappa$-H\"older  continuous for any $\kappa<H$. Its real
part $B:=\{ 2\Re\Gamma_t; \, t\in\R\}$  has the same law as fBm.

\item
The family of centered Gaussian real-valued processes $B(\eta):= \{ 2 \Re \Gamma_{t+\II\eta}; \, t\in\R\}$ converges a.s. to $B$ in $\alpha$-H\"older norm for any $\al<H$, on any interval $[0,T]$ with  $T>0$. Its infinitesimal
covariance kernel $\E [B'_x(\eta) B'_y(\eta)]$ is $
K'(\eta;x,y):=K^{',+}(\eta;x,y)+K^{',-}(\eta;x,y).
$

\end{enumerate}

\smallskip

\label{prop:Gamma-process}
\end{proposition}

\subsubsection{Definition of the L\'evy area}
Consider now an $m$-dimensional analytic fBm $\gga=(\gga^{(1)},\ldots,\gga^{(m)})$. Since the process $B(\eta)$ is  smooth, one can define the following integrals in the Riemann sense for all $0\le s<t \le T $, $1\le j_1,j_2\le m$ and $\eta>0$:
\begin{equation}\label{eq:8a}
\bb^{\bt,\eta}(j_1,j_2)
=\int_{s}^{t} dB_{u_1}^{(j_2)}(\eta) 
\int_{s}^{u_1} dB_{u_2}^{(j_1)}(\eta).
\end{equation}
It turns out that $\bb^{\bt,\eta}$ converges in the H\"older spaces $\cac_2^{2\ka}$ from Section \ref{sec:incr} (see \cite{Un,TU}),  which allows to define the L\'evy area in the following way:
\begin{proposition}\label{prop:2.3}
Let $T>0$ and define $\bb^{\bt,\eta}$ by equation (\ref{eq:8a}). Let also $0<\ga<H$. Then $B$ satisfies Hypothesis \ref{h1} in the following sense:
\begin{enumerate}
\item
The couple $(B(\eta),\bb^{\bt,\eta})$ converges in $L^p\big (\Omega;\cac_1^{\ga}(\ott;\R)\times\cac_2^{2\ga}(\ott^2;\R^{m,m}) \big)$ for all $p\ge 1$ to a couple $(B,\bb^\bt)$, where $B$ is a fractional Brownian motion.
 \item
The increment $\bb^\bt$ satisfies the  algebraic relation $\der\bb^\bt=\der B\otimes\der B$. 
\end{enumerate}
\end{proposition}

\medskip

One of the advantages of the analytic approach is that an expression for the covariances of the L\'evy area can  be easily derived by dominated convergence. We have
\begin{align}
\label{iso}
& \E  \lc \bb^{\bt}_{s_1t_1}(i,j) \,  \bb^{\bt}_{s_2t_2}(i,j)  \rc \\ \nonumber & \quad = H^2(2H-1)^2 \int_{s_1}^{t_1} \int_{s_2}^{t_2} \int_{s_1}^{u_1} \int_{s_2}^{u_2} |u_1-u_2 |^{2H-2} |v_1-v_2 |^{2H-2} \, dv_2 \, dv_1 \, du_2 \, du_1
\end{align}
for $0 \leq s_1 \leq t_1 \leq s_2 \leq t_2 \leq T$ and $i,j=1, \ldots, m$.

\smallskip

 Moreover,  $B(\eta)$  satisfies similar stationarity and scaling properties as the fBm itself.

\begin{lemma}\label{lemma_sc_st_B(eta)}
We have
\begin{enumerate}
\item (stationarity)
$$\lcl (\der B(\eta))_{s,u+s}, \ 0\leq u \leq T-s \rcl \stackrel{\mathcal{L}}{=} \lcl B(\eta)_u, \ 0\leq u\leq T-s \rcl, $$
\item (scaling) 
$$\lcl B(\eta)_{c \cdot u}, \ 0\leq u \leq T/c \rcl \stackrel{\mathcal{L}}{=} \lcl c^H B\lp \frac{\eta}{c}\rp_u, \ 0\leq u \leq T/c \rcl.$$
\end{enumerate}
\end{lemma}

The above Lemma can be shown by straightforward calculations exploiting that $B(\eta)$ is a Gaussian process with covariance kernel  $K'$ and
 will be useful to derive the 
   scaling property of the fractional L\'evy area. See Lemma \ref{stat} below.
 
\bigskip

\section{Approximation of the L\'evy area}

Let $ \mathcal{P}_{n,T}$ be the uniform partition $ \{ t_k^n=\frac{kT}{n}, k=0, \ldots, n \}$ of $\ott$, and let  $B^{n,T}$ be the linear interpolation of $B$ based on the points of $ \mathcal{P}_{n,T}$. More precisely, $B^{n,T}$ is defined as follows: for $t\in\ott$, let $k\in \{ 0, 1, \ldots , n-1\}$ be such that $t_k^n \leq t < t_{k+1}^n$. Then we have
\beq\label{eq:def-trapez-approx}
B^{n,T}_t=B_{t_k^{n}}+\lp \frac{t-t_k^n}{T/n}\rp (\der B)_{t_k^n t_{k+1}^n}.
\eeq 
Let also $\bb^{\bt,n,T}$ be the L\'evy area of $B^{n,T}$, which is simply defined in the Riemann sense by $$  \bb^{\bt,n,T}_{st}(i,j) =\ist \int_s^{u_1}dB^{n,T,(i)}_{u_2}\, dB^{n,T,(j)}_{u_1}. $$
The first step in the convergence analysis of our Milstein type scheme is to determine the  rate of convergence of the couple $(B^{n,T},\bb^{\bt,n,T})$ towards $(B,\bb^\bt)$. The current section is devoted to this step, which can be seen as an extension of \cite{NTU} to H\"older norms.  Throughout the remainder of this article we will denote unspecified non-negative and finite random variables by $\theta$, indicating by indices on which quantities they depend.
Similarly, we will denote unspecified constants, whose specific value is not relevant, by $C$ or $K$.

\subsection{Preliminary tools}

As a first preliminary step, let us state the following elementary lemma about the stationarity and scaling properties of the fBm $B$ and its piecewise linear interpolation $B^{n,T}$ resp. about the scaling property of the L\'evy areas $\bb^{\bt}$ and $\bb^{\bt,n,T}$.

\begin{lemma}\label{stat}
Consider a point  $s \in \mathcal{P}_{n,T}$. Then
\begin{equation}\label{stationarity_2}
\lcl (\der B)_{s,u+s},(\der B^{n,T})_{s,u+s}), \ 0\le u \leq T-s \rcl \stackrel{\mathcal{L}}{=} \lcl (B_u,B^{n,T}_u), \ 0\le u \leq T-s \rcl.
\end{equation}
Furthermore, if $c >0$, then
\begin{equation}\label{scaling_2}
\lcl (B_{cu},B^{n,T}_{cu}), \, 0\le u \leq T/c \rcl \stackrel{\mathcal{L}}{=} \lcl c^H(B_u,B_u^{n,T/c}), \ 0\le u \leq T/c \rcl.
\end{equation}
Finally, let $s,t \in \mathcal{P}_{n,T}$ with $s\leq t$. Then we have 
 \begin{align} \label{stat_scal_levy}
 \big( \bb^{\bt}_{st}(i,j) , \bb^{\bt,n,T}_{st}(i,j) \big)
  \, &\stackrel{\mathcal{L}}{=} \,  (t-s)^{2H}  \big(  \bb^{\bt}_{01}(i,j) , \bb^{\bt,n,T/(t-s)}_{01}(i,j) \big)
  \end{align}
  for all $i,j=1, \ldots, m.$ \end{lemma}

\begin{proof} 
These assertions are of course consequences of the stationarity and scaling properties of fBm, i.e.
for any $c>0$ the process
\begin{equation}\label{scaling}
 \widetilde{B}_{\cdot}^{(i)}=c^H B_{ \cdot /c}^{(i)}
\end{equation}
is again a fBm, and  for any $h\in\R$ the process
\begin{equation}\label{stationarity}
 \widetilde{B}_{\cdot}^{(i)}=(\der B^{(i)})_{h,\cdot+h}
\end{equation}
is a fBm. 

Recall that  the points of $\mathcal{P}_{n,T}$ are given by $t_i^n =\frac{iT}{n}$ for all $i\in \N$, and introduce the two mappings $F_-^{n,T}$ and $F_+^{n,T}$ defined on $\R_+$ by $F_-^{n,T}(u)=t_i^n$ and $F_+^{n,T}(u)=t_{i+1}^n$ if $t_i^n \leq u < t_{i+1}^n$. With these notations, one has $B^{n,T}_u=G^{n,T}(B)_u$, where the measurable mapping $G^{n,T}:  \cac( \R^+; \R^m) \rightarrow \cac( \R^+; \R^m)$ is defined by
$$G^{n,T}(y)_u=y_{F_-^{n,T}(u)}+\frac{u-F_-^{n,T}(u)}{T/n} \lp y_{F_+^{n,T}(u)}-y_{F^{n,T}_-(u)} \rp , \qquad u \in \R^+.$$
Now, in order to establish (\ref{stationarity_2}), note that $F_{\pm}^{n,T}(u+s)=F^{n,T}_{\pm}(u)+s$
if $s\in \mathcal{P}_{n,T}$. It is then easily seen that 
$$((\der B)_{s,\cdot +s},(\der B^{n,T})_{s,\cdot +s}) =((\der B)_{s,\cdot +s},G^{n,T}((\der B)_{s,\cdot +s}),$$
so that, due to the stationarity property of fBm, the following identity in law for processes holds true:
$$((\der B)_{s,\cdot +s},(\der B^{n,T})_{s,\cdot +s})\stackrel{\mathcal{L}}{=} (B,G^{n,T}(B)) =(B,B^{n,T}).$$

\smallskip

The proof of (\ref{scaling_2}) is quite similar. In fact, one has $F^{n,T}_{\pm}(c \cdot u)=c \cdot F^{n,T/c}_{\pm}(u)$ and so $B^{n,T}_{c u}=G^{n,T/c}(B_{c \cdot})_u$. Thus it holds, thanks to the scaling property of fBm,
$$(B_{c \cdot},B^{n,T}_{c \cdot})=(B_{c \cdot},G^{n,T/c}(B_{c \cdot})) \stackrel{\mathcal{L}}{=} (c^H B,G^{n,T/c}(c^H B)).$$
Identity (\ref{scaling_2}) is then a consequence of the linearity of $G^{n,T/c}$.

\smallskip

Now it remains to establish (\ref{stat_scal_levy}).
Note first that Proposition \ref{prop:2.3} implies
that
\begin{equation}  \label{conv_B_eta}
\big( { \bf B}^{\bf 2}_{st},  {\bf B}^{{\bf 2}, n,T }_{st}  \big) = 
\lim_{\eta \to 0} \big( { \bf B(\eta) }^{\bf 2}_{st}, { \bf B(\eta) }^{ {\bf 2}, n,T}_{st} \big) 
\end{equation}
 in probability.  Here ${\bf B(\eta)}_{st}^{ {\bf 2},n,T}$ is the L\'evy area associated to the piecewise linear interpolation of $B(\eta)$ with stepsize $T/n$.

Since $B(\eta)$ is analytic, the above L\'evy areas can be approximated by a standard Euler quadrature rule, i.e. we have
\begin{align}  \label{quad_B_eta_1}
{\bf B}(\eta)^{\bf 2}_{st}
= \lim_{ k \to \infty} \mathcal{I}_{k}({\bf B}(\eta)^{\bf 2}_{st}) \qquad 
{\bf B(\eta)}^{{\bf 2},n,T}_{st}
= \lim_{ k \to \infty} \mathcal{I}_{k}({\bf B(\eta)}^{{\bf 2},n,T}_{st})
\end{align}
almost surely, where
\begin{align*}
\mathcal{I}_{k}({\bf B}(\eta)^{\bf 2}_{st}) 
&= \sum_{i=0}^k \lcl (\der B(\eta))_{s,\frac{i}{k}(t-s)+s} \rcl \otimes \lcl (\der B(\eta))_{\frac{i}{k}(t-s)+s,\frac{i+1}{k}(t-s)+s}  \rcl \\
 \mathcal{I}_{k}({\bf B(\eta)}^{{\bf 2},n,T}_{st})
&=  \sum_{i=0}^k \lcl (\der B(\eta)^{n,T})_{s,\frac{i}{k}(t-s)+s} \rcl \otimes \lcl (\der B(\eta)^{n,T})_{\frac{i}{k}(t-s)+s,\frac{i+1}{k}(t-s)+s} \rcl . \end{align*}
Using again the $G^{n,T}$ notation and setting $\eta^{st}= \frac{\eta}{t-s}$
we have
\begin{align*}
& \big( \mathcal{I}_{k}({\bf B}(\eta)^{\bf 2}_{st}), \mathcal{I}_{k}({\bf B(\eta)}^{{\bf 2},n,T}_{st}) \big)= \\& \,\,  \Big( \sum_{i=0}^k \lcl (\der B(\eta))_{s,\frac{i}{k}(t-s)+s} \rcl \otimes \Big\{ \lc (\der B(\eta))_{s,\frac{i+1}{k}(t-s)+s} \rc  
-\lc (\der B(\eta))_{s,\frac{i}{k}(t-s)+s} \rc  \Big\},\\
&   
\sum_{i=0}^k G^{n,T}((\der B(\eta))_{s,.+s})_{\frac{i}{k}(t-s)} \otimes \Big\{ G^{n,T}((\der B(\eta))_{s,.+s})_{\frac{i+1}{k}(t-s)}
-G^{n,T}((\der B(\eta))_{s,.+s})_{\frac{i}{k}(t-s)} \Big\} \Big).
\end{align*}
Thus, invoking Lemma \ref{lemma_sc_st_B(eta)} and setting $\eta^{st}= \frac{\eta}{t-s}$, we end up with
\begin{align*}
& \big( \mathcal{I}_{k}({\bf B}(\eta)^{\bf 2}_{st}), \mathcal{I}_{k}({\bf B(\eta)}^{{\bf 2},n,T}_{st}) \big)\\
&\, \, 
\stackrel{\mathcal{L}}{=}  \Big( \sum_{i=0}^k  B(\eta)_{\frac{i}{k}(t-s)} \otimes   (\der B(\eta))_{\frac{i}{k}(t-s),\frac{i+1}{k}(t-s)}  ,\\
&  \hspace{5cm} \sum_{i=0}^k G^{n,T}(B(\eta))_{\frac{i}{k}(t-s)} \otimes  (\der G^{n,T}(B(\eta)))_{\frac{i}{k}(t-s),\frac{i+1}{k}(t-s)}  \Big) \\
& \, \, =  \Big( \sum_{i=0}^k  B(\eta)_{\frac{i}{k}(t-s)} \otimes  (\der  B(\eta))_{\frac{i}{k}(t-s),\frac{i+1}{k}(t-s)} ,\\
&  \hspace{5cm}  \sum_{i=0}^k G^{n,T/(t-s)}(B(\eta)_{.(t-s)})_{\frac{i}{k}} \otimes  (\der G^{n,T/(t-s)}(B(\eta)_{.(t-s)}))_{\frac{i}{k},\frac{i+1}{k}}  \Big) \\
&\,\, \stackrel{\mathcal{L}}{=} \Big((t-s)^{2H} \sum_{i=0}^k  B\lp \eta^{st} \rp_{\frac{i}{k}} \otimes   \lp \der B\lp \eta^{st} \rp\rp_{\frac{i}{k},\frac{i+1}{k}} ,\\
&  \hspace{3cm}   (t-s)^{2H} \sum_{i=0}^k G^{n,T/(t-s)}\lp B\lp \eta^{st}\rp\rp_{\frac{i}{k}} \otimes  \lp \der G^{n,T/(t-s)}\lp B\lp \eta^{st}\rp\rp\rp_{\frac{i}{k},\frac{i+1}{k}}  \Big)
,
\end{align*}
that is
\begin{align} \label{scal_app}
& \big( \mathcal{I}_{k}({\bf B}(\eta)^{\bf 2}_{st}), \mathcal{I}_{k}({\bf B(\eta)}^{{\bf 2},n,T}_{st}) \Big) \\ & \qquad \qquad \stackrel{\mathcal{L}}{=} (t-s)^{2H} \Big( \mathcal{I}_{k} \Big ({\bf B}\Big( \frac{\eta}{t-s} \Big)^{\bf 2}_{01} \Big), \mathcal{I}_{k} \Big ({\bf B}\Big( \frac{\eta}{t-s} \Big)^{{\bf 2},n,T/(t-s)}_{01} \Big ) \Big). \nonumber
\end{align}
Clearly, we also have
\begin{align}  \label{quad_B_eta_2}
{\bf B}\Big(  \frac{\eta}{t-s} \Big)^{\bf 2}_{01}
& = \lim_{ k \to \infty} \mathcal{I}_{k} \Big ( {\bf B}
\Big( \frac{\eta}{t-s} \Big)^{\bf 2}_{01} \Big), \\ \nonumber
{\bf B}\Big( \frac{\eta}{t-s} \Big)^{{\bf 2},n,T/(t-s)}_{01}
&= \lim_{ k\to \infty} \mathcal{I}_{k} \Big({\bf B}\Big( \frac{\eta}{t-s} \Big)^{{\bf 2},n,T/(t-s)}_{01} \Big)
\end{align}
 almost surely and
\begin{align}  \label{conv_B_eta_2}
\Big( { \bf B}^{\bf 2}_{01},  {\bf B}^{{\bf 2}, n,T/(t-s) }_{01}  \big) = 
\lim_{\eta \to 0} \Big( { \bf B } \Big( \frac{\eta}{t-s} \Big)^{\bf 2}_{01}, { \bf B}\Big(\frac{\eta}{t-s} \Big) ^{ {\bf 2}, n,T/(t-s)}_{01} \Big) 
\end{align} in probability.
So, combining  (\ref{conv_B_eta}), (\ref{quad_B_eta_1}), (\ref{scal_app}), (\ref{quad_B_eta_2}) and  (\ref{conv_B_eta_2}), we obtain
\begin{align*}
& \E \lc \vp \lp {\bf B}^{\bf 2}_{st}, {\bf B}^{ {\bf 2},n,T}_{st} \rp \rc  \\ & \qquad = \lim_{ k \rightarrow \infty, \eta \rightarrow 0} \E\lc \vp \big( \mathcal{I}_{k}({\bf B}(\eta)^{\bf 2}_{st}), \mathcal{I}_{k}({\bf B(\eta)}^{{\bf 2},n,T}_{st}) \big) \rc  \\ & \qquad  = 
\lim_{k \rightarrow \infty, \eta \rightarrow 0} \E \lc \vp \Big( (t-s)^{2H} \mathcal{I}_{k} \Big({\bf B}\lp \eta^{st} \rp^{\bf 2}_{01} \Big), (t-s)^{2H} \mathcal{I}_{k} \Big({\bf B}\lp \eta^{st} \rp^{{\bf 2},n,T/(t-s)}_{01} \Big) \Big) \rc \\ & \qquad = 
\E \lc \vp \lp (t-s)^{2H} {\bf B}^{ \bf 2}_{01}, (t-s)^{2H} {\bf B}^{ {\bf 2} ,n,T/(t-s)}_{01} \rp \rc
\end{align*}
for any function $\vp \in \cac_b((\R^m \otimes \R^m)^2)$, 
which concludes the proof of (\ref{stat_scal_levy}).

\end{proof}

\medskip

The next auxiliary result is an upper bound of the modulus of continuity of fBm and is a consequence of Theorem 3.1 in \cite{wang}.
 
 \begin{lemma}{\label{cont-fbm}}
 Let $T>0$. There exists  $h^*>0$ and a finite and non-negative random variable $\theta_{H,h^*,T}$ such that
 $$ \sup_{t \in \lc 0, T- h \rc } |(\der B)_{t,t+h}  | \leq \theta_{H,h^*,T} \cdot h^H \cdot \sqrt{|\log(1/h)|} $$
 for all $h \in (0,h^*).$
 \end{lemma}

\smallskip
 
The classical Garsia lemma reads as follows:
 \begin{lemma}\label{GRR-1} For all $\gamma >0$ and $ p \geq 1 $ there exists a constant $C_{\gamma,p,l}>0$  such that
\begin{align*} 
\cn \big [  f; \cac_1^{\gamma}([0,T]; \R^l) \big ]  \leq C_{\gamma,p,l} \left( \int_0^T \int_0^T  \frac{| (\delta f)_{uv}  |^{2p}}{|u-v|^{2  \gamma p+2}} \, du \,dv \right)^{1/(2p)}
\end{align*} 
for all $f \in \cac_1([0,T];\mathbb{R}^l)$.
\end{lemma}

\smallskip

Finally, we also need to control the H\"older smoothness
of elements of $\cac_2$, beyond the case of increments of functions in $\cac_1$. The following is a generalization  of the Garsia-Rodemich-Rumsey lemma above.

\begin{lemma}\label{GRR-2}Let $\kappa >0$ and $p \geq 1$. Let $R \in \cac_2([0,T]; \R^{l,l})$ with $\delta R \in \cac_3^{\ka}([0,T]; \R^{l,l})$.
If 
$$      \int_0^T \int_0^T  \frac{|R_{uv}|^{2p}}{|u-v|^{2 \kappa p +2}} \, du \,dv    < \infty, $$
then  $R \in \cac_2^{\kappa}([0,T]; \R^{l,l})$. In particular, there exists a constant $C_{\kappa, p,l}>0$,  such that
\begin{align*} 
& \cn \big [  R; \cac_2^{\kappa}([0,T]; \R^{l,l}) \big ]  \\ & \qquad  \leq C_{\kappa,p,l} \left( \int_0^T \int_0^T  \frac{| R_{uv}  |^{2p}}{|u-v|^{2  \kappa p+2}} \, du \,dv \right)^{1/(2p)} + C_{\kappa,p,l} \,\,  \cn \big[ \delta R; \cac_3^{\ka}([0,T]; \R^{l,l}) \big].
\end{align*} 
\end{lemma}

\smallskip
\medskip

\subsection{Approximation results}
Recall that our aim here is to show the convergence of the couple $(B^{n,T},\bb^{\bt,n,T})$ towards $(B,\bb^\bt)$ in some suitable H\"older spaces. A similar result was obtained in \cite{CQ}, but with the following differences: (i)~The authors in \cite{CQ} studied  the $p$-variation norm of  $\bb^\bt-\bb^{\bt,2^n,T}$ using dyadic discretisations, while we are working in the H\"older setting.  (ii) The rate of convergence for the approximation was not their main concern, and the convergence rate stated in \cite[Corollary 20]{CQ} is  not sharp.  

\smallskip

Let us now start with a first moment estimate for the difference $\bb^\bt-\bb^{\bt,n,T}$, for which we will use 
 the error bound for a trapezoidal approximation of $\bb^\bt$ derived in 
 \cite{NTU}. Moreover, recall that we denote  by $B^{n,T}$ the piecewise linear interpolation of $B$ on $[0,T]$ with respect to the uniform partition $ \mathcal{P}_{n,T}= \left \{t_{k}^{n};\, k=0, \ldots, n \right \}$, where $t_{k}^{n}= \frac{kT}{n}$, and by $\bb^{\bt,n,T}$ the corresponding L\'evy area.

\begin{proposition} \label{prop_NUT} Let $p \geq 1$ and $H>1/4$. Then, we have
$$ 
\left( \E \! \lln \bb^\bt_{0,T}-\bb^{\bt,n,T}_{0,T} \rrn^{p}  \right)^{1/p}
\leq K_p  \cdot T^{2H} \cdot n^{-2H+1/2}.$$
\end{proposition}

\begin{proof} 
First note that the random variable $\bb^\bt_{0,T}-\bb^{\bt,n,T}_{0,T}$ belongs to the sum of the first and the second chaos of $B$ (we refer to \cite{Nu-bk} for a specific description of these notions).
So all moments of  $\bb^\bt_{0,T}-\bb^{\bt,n,T}_{0,T}$ are equivalent and it suffices to show that there exists a constant $K>0$ such that, for all $i,j=1, \ldots, m$,
\beq\label{eq:l2-approx-B2}
\lp \E\lln \bb^\bt_{0,T}-\bb^{\bt,n,T}_{0,T} \rrn^{2} \rp^{1/2}
  \leq K \cdot T^{2H} \cdot n^{-2H+1/2}.
\eeq

\smallskip

Consider first the diagonal elements of $\bb^\bt_{0,T}-\bb^{\bt,n,T}_{0,T}$. In this case, we have $\bb^\bt_{0,T}(j,j)=(B_T^{(j)})^2/2$ and
\begin{align*}
\bb^{\bt,n,T}_{0,T}(j,j)&= \int_0^T B_u^{n,T,(j)} \, d B_u^{n,T,(j)}\\  
&= 
   \sum_{k=0}^{n-1}   B_{t_{k}^{n}}^{(j)} \der B_{t_{k}^{n}t_{k+1}^{n}}^{(j)}
   + \sum_{k=0}^{n-1} \int_{t_{k}^{n}}^{t_{k+1}^{n}}  \left( \frac{n}{T} \right)^2  \left (u - \frac{kT}{n}  \right)  \left( \der B_{t_{k}^{n}t_{k+1}^{n}}^{(j)}  \right)^2    \, du  \\
& = 
  \sum_{k=0}^{n-1}  \left(  B_{t_{k}^{n}}^{(j)} \, \der B_{t_{k}^{n}t_{k+1}^{n}}^{(j)}  
  + \frac{1}{2}  \left( \der B_{t_{k}^{n}t_{k+1}^{n}}^{(j)}  \right)^2  \right) 
 = \frac{1}{2} \left( B_T^{(j)}  \right)^2.
 \end{align*}
Hence it follows
\beq\label{eq:approx-B2-diag}
\bb^\bt_{0,T}(j,j)-\bb^{\bt,n,T}_{0,T}(j,j)=\int_0^T B_u^{(j)} \, d B_u^{(j)} - \int_0^T B_u^{n,T,(j)} \, d B_u^{n,T,(j)} =0.
\eeq

\smallskip

Now consider the off-diagonal terms  of $\bb^\bt_{0,T}-\bb^{\bt,n,T}_{0,T}$. Without loss of generality we can assume that $i>j$. Proceeding as above we have
 $$ \int_0^T B_u^{n,T,(i)} \, d B_u^{n,T,(j)} = \frac{1}{2} \sum_{k=0}^{n-1}  
 \big( B_{t_{k}^{n}}^{(i)}  + B_{t_{k+1}^{n}}^{(i)} \big) \,  
 \der B_{t_{k}^{n}t_{k+1}^{n}}^{(j)}.
 $$
Thus, 
 \cite[Theorem 1.2]{NTU} can be applied and yields
\begin{align}
\left( \E \left|   \bb^\bt_{0,T}(i,j)   -   \bb^{\bt,n,T}_{0,T}(i,j) \right |^2 \right)^{1/2}  \leq \, K  \cdot T^{2H} \cdot n^{-2H+1/2}.
\end{align}

\end{proof}

\smallskip
 The next result gives an  error bound for the piecewise linear interpolation of $B$. Note that similar estimates as in the next lemma can be found in \cite{DN}, where the case $H>1/2$ is considered. 
 
 \begin{lemma}\label{lemma_lin_int}
Let $0 \leq   \gamma < H. $ Then, there exists a finite and non-negative random variable
 $\theta_{H, \gamma, T}$ such that 
$$  \cn \big[  B^{n,T}-B \, ;\cac_1^\ga([0,T]) \big] \leq \theta_{H, \gamma, T} \cdot \sqrt{\log(n)} \cdot n^{-(H-\gamma)}$$ for $n>1$.
\end{lemma}
 
 \begin{proof}
 Clearly, we have to find appropriate bounds for 
 $$ |\der (B^{n,T}-B)_{st}|, \qquad s,t \in [0,T].$$
First note that there exists a strictly positive $x_{H, \ga}$ such that the mapping $f:(0,T] \rightarrow [0, \infty)$, $f(x)=x^{H-\ga} \sqrt{ |\log(1/x) |}$ is increasing on $x \in (0, x_{H, \ga})$. Without loss of generality, we assume that $T/n \leq \inf(x_{H, \ga},h^*)$, where $h^*$ is defined by Lemma \ref{cont-fbm}. 

 \smallskip
 
 \noindent
(i) First, consider the case where $\lln t-s \rrn \geq \frac{T}{n}$. Let us assume also without loss of generality that $t_k^n \leq s < t_{k+1}^n \leq t_l^n \leq t < t_{l+1}^n$ for some $k<l$ and recall that $t_k^n=kT/n$. Then 
$$
B^{n,T}_s=B_{t_k^{n}}+\lp \frac{s-t_k^n}{T/n}\rp \der B_{t_k^n t_{k+1}^n}
\quad\mbox{and}\quad
B^{n,T}_t=B_{t_l^n}+\lp \frac{t-t_l^n}{T/n}\rp \der B_{t_l^n t_{l+1}^n},
$$
so that
\begin{align*}
 | \der(B^{n,T}-B)_{st}| & \leq     |\der B_{t_k^n s}| + |\der B_{t_n^l t}| + | \der B_{t_k^n t_{k+1}^n}|  +  |\der B_{t_l^n t_{l+1}^n}|  \\  &   \leq  4 \theta_{H, T} \sqrt{ |\log(n/T) |} 
  \left( \frac{ T}{n} \right)^{H} \leq \theta_{H, T} |t-s|^{\gamma} \sqrt{ |\log(n)|} n^{-(H-\ga)} 
\end{align*}
using Lemma  {\ref{cont-fbm}}.

 \smallskip
 
 \noindent
(ii) Now, suppose that $\lln t-s \rrn < T/n$ with for instance $t_k^n \leq s < t < t_{k+1}^n$. In this case, $$(\der B^{n,T})_{st}=\frac{t-s}{T/n}(\der B)_{t_k^n t_{k+1}^n}$$ and thus
\begin{align*}   & |  \der(B^{n,T}-B)_{st}|  
\leq |\der B_{st}| +   |  \der B^{n,T}_{st}|
\\  
& \qquad \qquad \leq \theta_{H, T}   \sqrt{ |\log(1/(t-s)) |} |t-s|^{H} + \theta_{H, T}  |t-s| \sqrt{|\log(n/T)|}\left( \frac{T}{n} \right)^{H-1}  \\
&  \qquad \qquad  \leq  \theta_{H, T}   \sqrt{ |\log(1/(t-s)) |} |t-s|^{H} + \theta_{H, T}  |t-s|^{\ga} \sqrt{|\log(n)|}n^{-(H-\ga)}  .\end{align*}
Using the monotonicity of $x \mapsto x^{H-\ga} \sqrt{ |\log(1/x) |}$,
it follows
\begin{align*}
 | \der(B^{n,T}-B)_{st} | \leq {\theta}_{H,T} |t-s|^{\gamma} \sqrt{|\log(n) |} n^{-(H-\ga)}. 
\end{align*}

 \smallskip
 
 \noindent
(iii) The same estimate as above also holds true if $\lln t-s\rrn < T/n$ and $t_k^n \leq s < t_{k+1}^n \leq t < t_{k+2}^n$.

 \smallskip
 
 \noindent
(iv) Combining (i)-(iii) yields the assertion.

 \end{proof}

 Now we determine the error for the approximation of the L\'evy area.
 
\begin{lemma}\label{main_wz_Lp}
Let $1/4 <  \gamma < H$. Then, there exists a finite and non-negative random variable
 $\theta_{H,  \gamma, T}$ such that 
$$  \cn \big[ {\bf B}^{ {\bf 2},n,T}-\bd \, ;\cac_2^{2\ga}([0,T])\big] 
\leq
 \theta_{ H, \gamma,  T} \cdot \sqrt{\log(n)} \cdot n^{-(H-\gamma)} $$ for $n>1$.
\end{lemma}

\begin{proof} In this proof we will denote constants (which depend only on  $p$, $q$, $\varepsilon$, $\gamma $  and $T$) by $K$, regardless of their value. 
 
 \smallskip
 
 \noindent
{\it Step 1.}  We will first show that 
\begin{align}\label{result-x-2}
\left( \E \left|\cn[ {\bf B}^{ {\bf 2},n,T}-\bd;\cac_2^{2\ga}([0,T])] \right|^q\right)^{1/q} \leq K \cdot \big( n^{-2(H-\gamma)} + n^{-H} \big).
\end{align}

For this, we have to consider the family of increments $A^{n,T}(i,j)\in\cac_2$, defined by
$$ A^{n,T}_{st}(i,j)= \int_s^t (\der B^{(i)})_{su}  \, d B_u^{(j)}   - \int_s^t (\der B^{n,T,(i)})_{su}   \, d B_u^{n,T,(j)} $$
for $i,j=1, \ldots, m$. By symmetry we can assume $1 \leq j\leq i \leq m$.

 \smallskip

We distinguish several cases for $s,t \in \lc 0, T \rc$.

\smallskip
 
 \noindent
(i) Assume that $\lln t-s \rrn \geq \frac{T}{n}$ and $s,t \in \mathcal{P}_{n,T}$, i.e. $s=\frac{kT}{n}$ and  $t=\frac{lT}{n}$ for $k<l$. Then the scaling properties of fBm, see Lemma \ref{stat}, yield 
\begin{align*}
 A^{n,T}_{st}(i,j)  & \stackrel{\mathcal{L}}{=} \int_0^{t-s} B_u^{(i)} \, d B_u^{(j)}   - \int_0^{t-s} B_u^{n,T,(i)}   \, d B_u^{n,T,(j)} \\ 
 & \stackrel{\mathcal{L}}{=}  (t-s)^{2H} \lp \int_0^1 B_u^{(i)}  d B_u^{(j)} -\int_0^1 B^{n,T/(t-s),(i)}_u \,d  B^{n,T/(t-s),(j)}_u \rp.
\end{align*}
Since $\frac{T}{t-s}=\frac{n}{l-k}$
we have
$$ \left \{ B^{n,T/(t-s),(i)}_u, \,  u \in [0,1] \right \} = \left \{ B^{l-k,1,(i)}_u, \, u \in [0,1] \right \}.$$

Now Proposition \ref{prop_NUT} gives 
\begin{align} \label{levy_disc_error} \nonumber
\left( \E \left| A^{n,T}_{st}(i,j) \right|^p \right)^{1/p} & \leq K \cdot |t-s|^{2H} \cdot |l-k|^{-2H+1/2} \leq K \cdot |t-s|^{1/2} \cdot n^{-2H+1/2} \\ 
& \leq K \cdot |t-s|^{2\gamma}  \cdot n^{-2(H-\gamma)} ,
\end{align}
with $\ga>1/4$. 

 \smallskip
 
 \noindent
(ii) Assume now that $(t-s) \geq \frac{T}{n}$ with $ s < t_{k+1}^n \leq t_l^n \leq t < t_{l+1}^n$. 
Using  the cohomologic relation $\der(\der A^{n,T}(i,j))_{s t_{k+1}^n t_l^n t}=0$, we obtain
\begin{multline}\label{decompo-suggest}
A^{n,T}_{st}(i,j)= A^{n,T}_{st_{k+1}^n}(i,j)+A^{n,T}_{t_{k+1}^n t_l^n}(i,j)+A^{n,T}_{t_l^n t}(i,j)\\
 +\der(A^{n,T}(i,j))_{st_{k+1}^n t}+\der(A^{n,T}(i,j))_{t_{k+1}^nt_l^n  t} .
\end{multline}

For the term $A^{n,T}_{t_{k+1}^n t_l^n}(i,j)$, we can use the first step to deduce
$$ \left ( \E  \left |A^{n,T}_{t_{k+1}^n t_l^n}(i,j) \right|^p \right) ^{1/p} \leq K \cdot \lln t_l^n -t_{k+1}^n \rrn^{2\ga}\cdot n^{-2(H-\ga)} \leq K \cdot \lln t-s \rrn^{2\ga} \cdot n^{-2(H-\ga)}.$$
To deal with the last two terms of (\ref{decompo-suggest}), remember the algebraic relation
\begin{equation}\label{rel-alg-proof}
\der(A^{n,T}(i,j))=\der B^{(i)} \cdot \der B^{(j)}-\der B^{n,T,(i)} \cdot \der B^{n,T,(j)},
\end{equation}
which entails here
$$| \der(A^{n,T}(i,j))_{st_{k+1}^n t}| \leq |(\der B^{(i)})_{st_{k+1}^n}| \cdot |(\der B^{(j)})_{t_{k+1}^n t } | +|(\der B^{n,T,(i)})_{st_{k+1}^n}| \cdot |(\der B^{n,T,(j)})_{t_{k+1}^n t} |,$$
and we easily get
$$ \left ( \E \left |\der (A^{n,T}(i,j))_{st_{k+1}^n t} \right|^p \right)^{1/p}\leq K \cdot \lln t-s \rrn^H \cdot (T/n)^H  \leq K \cdot \lln t-s \rrn^{2\ga} \cdot \big( n^{-2(H-\ga)} + n^{-H} \big).$$

Similarly we obtain the same estimate for $\E [ |\der (A^{n,T}(i,j))_{t_{k+1}^n t_l^n t}|^p ]^{1/p}$.

\smallskip

As for the term $A^{n,T}_{st_{k+1}^n}(i,j)$  one has, on the one hand,
\begin{align}\label{eq:44}  \left(  \E \left|  \int_s^{t_{k+1}^n}(\der B^{(i)})_{su} \,  dB^{(j)}_u \right |^p \right)^{1/p} &= \lln t_{k+1}^n -s \rrn^{2H} \left(  \E \left|  \int_0^{1}B^{(i)}_u\, dB^{(j)}_u \right |^p \right)^{1/p}  \\ & \leq K  \cdot \lln t -s \rrn^{2\gamma} \cdot n^{-2(H-\gamma)},
\end{align}
where $\ga<H$.
On the other hand, 
$$
\left| \int_s^{t_{k+1}^n} \der B^{n,T,(i)}_{su}\, dB^{n,T,(j)}_u \right|
=\left| \der B^{n,T,(i)}_{t_k^n t_{k+1}^n} \, \der B^{n,T,(j)}_{t_k^n t_{k+1}^n} \right| 
\int_s^{t_{k+1}^n} \frac{(u-t_k^n)}{(T/n)^2} du
\le \left|  \der B^{n,T,(i)}_{t_k^n t_{k+1}^n} \, \der B^{n,T,(j)}_{t_k^n t_{k+1}^n} \right|.
$$ 
So for $\ga<H$, an application of the Cauchy-Schwarz inequality yields
\begin{equation}\label{eq:45}  \left( \E \left|  \int_s^{t_{k+1}^n}(\der B^{n,T,(i)})_{su} \,  dB^{n,T,(j)}_u \right |^p \right)^{1/p}  \leq  K  \cdot \lln t -s \rrn^{2\gamma} \cdot n^{-2(H-\gamma)}. 
\end{equation}
Putting together relation (\ref{eq:44}) and  (\ref{eq:45}), we obtain $ (\E  [|A^{n,T}_{st_{k+1}^n }(i,j) |^p])^{1/p}   \leq K \lln t-s \rrn^{2\ga} \cdot n^{-2(H-\ga)}$. Furthermore, the term $A^{n,T}_{t_{l}^nt}(i,j)$ can be handled along the same lines.

 \noindent
(iii) It only remains to analyze the case $(t-s) < \frac{T}{n}$. For $t_k^n \leq s < t < t_{k+1}^n$ we have\begin{align*}
 \left( \E\left|\int_s^{t}(\der B^{(i)})_{su}\,  dB^{(j)}_u \right|^p \right)^{1/p} \leq  K \cdot | t -s |^{2H}\leq K \cdot |t-s|^{2\gamma}\cdot n^{-2(H-\gamma)},
\end{align*}
and
$$
\left | \int_s^t \der B^{n,T,(i)}_{su} \, dB^{n,T,(j)}_u \right| 
\le 
\frac{(t-s)^2}{2 (T/n)^2} \, 
\left| \der B_{t_k^n t_{k+1}^n}^{n,T,(i)} \right|  \,  \left|
\der B_{t_k^n t_{k+1}^n}^{n,T,(j)} \right|,
$$
and thus
$$
 \left( \E \left|\int_s^t (\der B^{n,T,(i)})_{su} \, dB^{n,T,(j)}_u \right|^p \right)^{1/p} \leq K \cdot |t-s|^{2 \gamma}\cdot n^{-2(H-\gamma)}.
$$
The case $(t-s) < \frac{T}{n}$ and $t_k^n \leq s < t_{k+1}^n \leq t < t_{k+2}^n$ can be treated analogously.

 \smallskip
 
 \noindent
(iv) Combining  steps (i)--(iii)  yields that
\begin{align} \label{levy_error_total}
\left( \E \left| A^{n,T}_{st}(i,j) \right|^p \right)^{1/p}  
 \leq K \cdot |t-s|^{2\gamma}  \cdot \big( n^{-2(H-\gamma)} + n^{-H} \big)  
 \end{align}
for all $s,t \in [0,T]$ and $1/4<\gamma <H$.

\bigskip

 \noindent
{\it  Step 2.} Before we can apply Lemma \ref{GRR-2}, we need additional preparations. First, notice that (\ref{rel-alg-proof}) can also be written as
$$\der(\bd-{\bf B}^{ {\bf 2},n,T})=\lc \der \lp B-B^{n,T}\rp \rc \otimes \der B + \der B^{n,T} \otimes \lc \der \lp B-B^{n,T} \rp \rc,
$$
so that 
\begin{multline*}
|\der(\bd-{\bf B}^{ {\bf 2},n,T})_{sut}| \\
\leq |  t-u|^{\ga}  |s-u|^{\ga} \left( 2\cn[\der B;\cac_2^\ga] \cdot  \cn[\der (B-B^{n,T});\cac_2^\ga]+ ( \cn[\der (B-B^{n,T});\cac_2^\ga])^2\right)
\end{multline*}
and thus
\begin{align*} 
\cn \big[  \der(\bd-{\bf B}^{ {\bf 2},n,T}); \cac_3^{2\gamma} \big] \leq  2 \cn[\der B;\cac_2^\ga] \cdot  \cn[\der (B-B^{n,T});\cac_2^\ga]+ ( \cn[\der (B-B^{n,T});\cac_2^\ga])^2.
\end{align*}

Lemma \ref{lemma_lin_int} now gives
\begin{align} \label{GRR-2-app-1} 
 \cn \big[  \der(\bd-{\bf B}^{ {\bf 2},n,T}) ; \cac_3^{2\ga} \big]  \leq  \theta_{H,\ga, T} \cdot  \sqrt{\log(n)} \cdot n^{-(H-\gamma)}.
\end{align}
\smallskip
 
 \noindent
{\it  Step 3.} Using (\ref{GRR-2-app-1}), Lemma \ref{GRR-2} entails

\begin{align*} 
& \cn \big [ (\bd-{\bf B}^{ {\bf 2},n,T}) ; \cac_2^{2\gamma}([0,T] \big ]  \\ & \qquad  \leq K \left( \int_0 ^T \int_0^T  \frac{| (\bd-{\bf B}^{ {\bf 2},n,T})_{uv}  |^{2p}}{|u-v|^{4 \gamma p+2}} \, du \,dv \right)^{1/(2p)} + K \cdot \theta_{H, \ga, T}\cdot  \sqrt{\log(n)} \cdot n^{-(H-\gamma)}.
\end{align*} 
for all $p\geq 1$.
To finish the proof, it remains to show that
\begin{align}  \label{to_show_bc}
 |R_{n,p}| \leq  \theta_{\ga,H,T}\cdot  \sqrt{\log(n)} \cdot n^{-(H-\gamma)} \end{align}
 where
 $$ R_{n,p}= \left( \int_0 ^T \int_0^T  \frac{| (\bd-{\bf B}^{ {\bf 2},n,T})_{uv}  |^{2p}}{|u-v|^{4 \gamma p+2}} \, du \,dv \right)^{1/(2p)}.$$
However, using (\ref{levy_error_total}) with $\gamma+ \varepsilon/2$ instead of $\gamma$, we have
\begin{align*}
 \E  |R_{n,p}|^{2p} & \leq    \int_0 ^T \int_0^T   \frac{ \E | (\bd-{\bf B}^{ {\bf 2},n,T})_{uv}  |^{2p}}{|u-v|^{4 \gamma p+2}} \, du \,dv  \\
 & \leq K    \int_0 ^T \int_0^T   \frac{| u-v|^{4 \gamma p+ 2 \varepsilon p}}{|u-v|^{4 \gamma p+2}} \, du \,dv \cdot \big( n^{-4(H-\gamma-\varepsilon/2)p} + n^{-2Hp} \big),
\end{align*}
i.e.
\begin{align*}
 (\E  |R_{n,p}|^{2p})^{1/(2p)} 
  \leq K    \int_0 ^T \int_0^T   | u-v|^{ 2p \varepsilon -2 } \, du \,dv \cdot \big( n^{-2(H-\gamma)+\varepsilon} + n^{-H} \big).
\end{align*}
So for $ p > \frac{1}{ \varepsilon}$, it holds
\begin{align*}
 (\E  |R_{n,p}|^{2p})^{1/(2p)} 
  \leq K    \cdot \big( n^{-2(H-\gamma)+\varepsilon} + n^{-H} \big).
\end{align*}
Now, set $ \alpha=\min\{2(H-\gamma)-\varepsilon, H \}$ and let 
$\delta>0$. From the
Chebyshev-Markov inequality  it follows
$$ {\bf P}( n^{\alpha-\varepsilon} |R_{n,p}| > \delta) \leq \frac{\E
|R_{n,p}|^{2p}}{\delta^{2p}}  n^{ 2p(\alpha - \varepsilon)} \leq K \frac{n^{-2p\varepsilon}}{\delta^{2p}} .$$
Since $p>1/\varepsilon$ we have
$$ \sum_{n=1}^{\infty} {\bf P}(n^{\alpha-\varepsilon} |R_{n,p}| > \delta) < \infty $$
for all $\delta >0$. The Borel-Cantelli Lemma   implies  now that $ n^{ \alpha - \varepsilon} |R_{n,p}|
\rightarrow 0$ a.s. for $n \rightarrow  \infty$, which gives   (\ref{to_show_bc}) by choosing $\varepsilon >0 $ appropriately, since
$$ \alpha-\varepsilon= \min \{ 2(H-\gamma - \varepsilon), H - \varepsilon \} > H- \gamma.$$

\end{proof}

\bigskip

Recall that the  Wong-Zakai approximation $\overline{Z}^n$ of $Y$ has been defined at equation (\ref{WZ}) by
\beq \label{wz_2}
\overline{Z}^n_t= a+ \sum_{i=1}^m\int_0^t \sigma^{(i)}(\overline{Z}_u^n) \, d B_u^{(i),n,T}, \quad t \in \lc 0, T \rc, \qquad  a \in \R^d.
\eeq
In particular, $\overline{Z}^n$ can be expressed as
$\overline{Z}^n =F(a,  B^{n,T} ,{\bf B}^{ {\bf 2},n,T})$, using Theorem~\ref{thm:Lip}. Hence, as a direct application of Lemmata \ref{lemma_lin_int} and \ref{main_wz_Lp} and invoking the Lipschitzness of $F$, we obtain  the following error bound for the Wong-Zakai approximation.

\begin{proposition} \label{step_1} Let $T>0$ and $1/3 <\gamma < H$. Then, there exists a finite random variable $\eta_{H ,\gamma, \si, T}^{(1)}$ such that
$$      \| Y- \overline{Z}^{n} \|_{\ga, \infty,T}   \leq  \eta_{H,\gamma, \si, T}^{(1)} \cdot \sqrt{\log(n)} \cdot n^{-(H-\gamma)}$$ for $n>1$.
\end{proposition}

\bigskip

\section{Discretising the Wong-Zakai approximation}

In the last section we have  established an error bound for  the Wong-Zakai approximation $\overline{Z}^{n}$ of the real solution $Y$. As mentioned in the introduction, the Milstein scheme corresponding to $\overline{Z}^{n}$ is exactly our simplified Milstein scheme (\ref{eq:misltein-scheme_2}). Thus, it remains to determine the discretisation error for $\overline{Z}^{n}$ itself.  To this aim, we first give  a general error bound for the Milstein scheme for ordinary differential equations (ODEs) driven by a smooth path $x$. Since Theorem 
\ref{thm:Lip} allows to  derive a non-classical stability result (in $\gamma$-H\"older norm) for the flow of an ODE driven by a smooth path, we can follow here the techniques of the numerical analysis for classical ODEs. In a second step, we will apply these bounds to our particular fBm approximation.

\subsection{The Milstein scheme for ODEs driven by smooth paths}
In this section, consider a piecewise  differentiable path $x \in C (\lc 0 , T \rc ;  \R^l)$ and  a function $g \in C^{3}(\mathbb{R}^d; \mathbb{R}^{d,l})$ which is bounded with bounded derivatives.
For the ordinary differential equation
 \begin{align} \label{ODE}
 \dot{y}_t= \sum_{i=1}^l g^{(i)}(y_t) \, dx^{(i)}_t , \quad t \in \lc 0, T \rc, \qquad a \in \mathbb{R}^d,
 \end{align}
the classical second order Taylor scheme with stepsize $T/n$ reads as:  $z^{n}_{0}=a$ and
\begin{equation}\label{eq:taylor-scheme}
z^{n}_{k+1} =  z^{n}_{k}   
+   \sum_{i=1}^lg^{(i)}(z^{n}_{k}) \, \der x^{(i)}_{t_k t_{k+1}} 
+   \sum_{i,j=1}^l   \mathcal{D}^{(i)} g^{(j)}(z^{n}_{k}) \int_{t_k}^{t_{k+1}} \, \der x^{(i)}_{t_k s} \, d x^{(j)}_s,
\end{equation}
where $\mathcal{D}^{(i)} = \sum_{p=1}^d  g_p^{(i)} {\partial}_p$, and where we have set $z^{n}_{k}=z^{n}_{t_k^n}$ with $t_k^n=kT/n$. For notational simplicity we will write in the following $t_k$ instead of $t_k^n$.
Introducing the numerical flow
\begin{equation}\label{eq:def-num-flow}
\Psi(z;t_k,t_{k+1}):= z 
+   \sum_{i=1}^lg^{(i)}(z) \, \der x^{(i)}_{t_k t_{k+1}}  
+   \sum_{i,j=1}^l   \mathcal{D}^{(i)} g^{(j)}(z) \int_{t_k}^{t_{k+1}} \, \der x^{(i)}_{t_k s} \, d x^{(j)}_s
\end{equation}
we can write this scheme as
\begin{align*} 
 z^n_{0}=a,  & \qquad
z^n_{t_{k+1}} =  \Psi(z^n_{k};t_{k},t_{k+1}) , \quad k=0, \ldots, n-1  .
\end{align*}
 For $q>k$ we also define
$$ \Psi(z;t_k,t_q):= \Psi( \cdot ;t_{q-1},t_q) \circ  \Psi( \cdot ;t_{q-2},t_{q-1}) \circ \cdots \circ \Psi( z ; t_k,t_{k+1}) .$$

 Moreover, the flow 
$\Phi(z;s,t)$ of the ODE (\ref{ODE}) is given by $\Phi(z;s,t):=y_t$,  where $y$ is the unique solution of 
\begin{align} \label{flow}
 \dot{y}_t= \sum_{i=1}^l g^{(i)}(y_t) \, dx^{(i)}_t , \quad t \in \lc s, T \rc, \qquad y_s=z.
 \end{align}

\smallskip
 
A straightforward Taylor expansion of the flow of the ODE gives that the one-step error
$$ r_k = \Phi(z;t_k,t_{k+1})- \Psi(z;t_k,t_{k+1}) $$
satisfies
\begin{align} \label{one-step-bound}
| r_k| &\leq  C \cdot  \sup_{i,j,p=0, \ldots, m} \|    \mathcal{D}^{(i)}   \mathcal{D}^{(j)} g^{(p)}\|_{\infty} \cdot  M_{t_kt_{k+1}}^{x}
\end{align}
with
\begin{align*}
M_{st}^{x}:=\left |\int_{s}^{t} |\dot{x}_w| \, dw \right|^3. 
\end{align*}

\smallskip

Furthermore, considering the smooth path $x$ as a rough path, Theorem \ref{thm:Lip} directly yields the following stability result for the flow:
 
\begin{proposition}\label{davie_lemma-bis} Let  $1/3 <  \ga  \leq 1$ 
and set $\| {\bf x }\|_{\ga} = \|  x \|_{\ga} + \| {\bf x^2} \|_{2\ga}$.
Then, there exists an increasing function  $C_T: \R \rightarrow \R_+$  such that
\begin{align}  \label{stab2-bis}  \frac{\left | (\Phi(z;s,t) -\Phi(\tilde{z};s,t))-(z- \tilde{z}) \right|}{ |t-s|^{\ga } } \leq  
C_T(\| {\bf x }\|_{\ga}) \cdot |z - \tilde{z}| \end{align}
and
\begin{align} \label{stab1-bis} \left |\Phi(z;s,t) -\Phi(\tilde{z};s,t) \right| \leq  C_T(\| {\bf x }\|_{\ga}) \cdot |z - \tilde{z}|  \end{align}
for all $s,t \in [0,T] $ and $z, \tilde{z} \in \mathbb{R}^d$.
\end{proposition}

\smallskip

The following stability result is crucial to derive the announced error bound for the Milstein scheme.

\begin{proposition}\label{prop:apriori}
Let $x \in C (\lc 0 , T \rc ;  \R^l)$ be a piecewise  differentiable path, and $g \in C_b^{3}(\mathbb{R}^d; \mathbb{R}^{d,l})$. Consider the flow $\Phi$ given by equation (\ref{flow}) and the numerical flow $\Psi$ defined by relation (\ref{eq:def-num-flow}). For $k=0,\ldots,n$, let $t_k=kT/n$, $y_{t_k}=\Phi(a;0,t_k)$ and $z_{t_k}=\Psi(a;0,t_k)$. Moreover recall that we have set
$$ M_{st}^{x}= \left|\int_s^t |\dot x_u| \, du \right|^{3}, \qquad 0\le s<t\le T.$$ Then, there exists an increasing function $\tilde{C}_T: \R \rightarrow \R^+$ such that  we have
\begin{eqnarray}
|y_{t_q} - z_{q}^n| &\leq&   \tilde{C}_T(\| {\bf x }\|_{\ga}) \cdot \sum_{k=0}^{q-1} 
M_{t_k t_{k+1}}^{x}
 \label{ode-error-bis} \\
|\der (y-z^n)_{t_p t_q}  |    &\leq  &
\tilde{C}_T(\| {\bf x }\|_{\ga}) \cdot \lcl  \sum_{k=p}^{q-1} M_{t_k t_{k+1}}^{x} +   
|t_q-t_p|^{\ga}  \cdot \sum_{k=0}^{p-1} M_{t_k t_{k+1}}^{x} \rcl
\label{ode-error-2-bis}
\end{eqnarray}
for $0\le p<q\le n$.

\end{proposition}

\begin{proof}
We will use the classical decomposition of the error in terms of the exact and the numerical flow: 
Since
$ z_{k}^{n}= \Phi(z_{k}^{n}; t_{k},t_{k})$ and $y_{t_k}=\Phi(z_{0}^{n};t_0,t_k)$, one has
$$y_{t_q} - z_{q}^{n} = \Phi(z_{0}^{n};t_0,t_q)-\Phi(z_{q}^{n}; t_{q},t_{q}) =
\sum_{k=0}^{q-1} \big( \Phi(z_{k}^{n};t_k,t_q)  - \Phi(z_{k+1}^{n};t_{k+1},t_q) \big). $$
Furthermore, thanks to the relation
$$ \Phi(z_{k}^{n};t_{k},t_q) = \Phi( \Phi(z_{k}^{n};t_k,t_{k+1});t_{k+1},t_q),$$
 the stability result (\ref{stab1-bis})  implies
$$  \big| \Phi(z_{k}^{n};t_k,t_q)  - \Phi(z_{k+1}^{n};t_{k+1},t_q) \big|          \leq 
{C}_T(\| {\bf x }\|_{\ga}) \cdot \left| \Phi(z_{k}^{n};t_k,t_{k+1}) - z_{{k+1}}^{n} \right|.
$$
However,  (\ref{one-step-bound}) gives
$$   \left| \Phi(z_{k}^{n};t_k,t_{k+1}) - z_{k+1}^{n} \right| = \left| \Phi(z_{k}^{n};t_{k},t_{k+1}) - \Psi(z_{k}^n;t_k,t_{k+1}) \right|  \leq C \cdot M_{t_k t_{k+1}}^{x},
 $$
from which (\ref{ode-error-bis}) is easily deduced.

\smallskip

Moreover, for $q \geq p$ we also have
\begin{align*}
&\der (y-z^n)_{t_p t_q}  = \lp\Phi(y_{t_p};t_p,t_q) -y_{t_p}\rp-\lp\Psi(z^n_{p};t_p,t_q) - z^n_{p}\rp \\
& =   \lp\Phi(y_{t_p};t_p,t_q) -y_{t_p}\rp - \lp\Phi(z^n_{p};t_p,t_q) -z^n_{p})\rp   - \lp\Psi(z^n_{p};t_p,t_q) - \Phi(z^n_{p};t_p,t_q)\rp.
\end{align*}
Analogously to  the derivation of (\ref{ode-error-bis}), one can show that
\begin{align} \label{ode-error-2h-bis}
| \Psi(z^n_{p};t_p,t_q) - \Phi(z^n_{p};t_p,t_q)| \leq  C \cdot {C}_T(\| {\bf x }\|_{\ga}) \cdot \sum_{k=p}^{q-1} M_{t_k t_{k+1}}^{x}.
\end{align}
Using (\ref{stab2-bis}) and (\ref{ode-error-bis}) we trivially end up with (\ref{ode-error-2-bis}).

\end{proof}

\medskip

\subsection{Application to fBm}
In order to apply Proposition \ref{prop:apriori} to the Wong-Zakai approximation    $\overline{Z}^n$  given by (\ref{wz_2}) note once again that our Milstein-type scheme $Z_{t_{0}}^{n}=a$ and
 \begin{equation*} 
Z_{t_{k+1}}^{n} =  Z_{t_k}^{n}   
+  \sum_{i=1}^m \sigma^{(i)} ( Z_{t_{k}}^{n}) \, \der B^{(i)}_{t_{k}t_{k+1}}
+  \frac{1}{2} \sum_{i,j=1}^m \mathcal{D}^{(i)}\sigma^{(j)} ( Z_{t_{k}}^{n})
\, \der B^{(i)}_{t_{k}t_{k+1}} \, \der B^{(j)}_{t_{k}t_{k+1}}
\end{equation*}
is obtained by discretising the Wong-Zakai approximation with the standard second order Taylor scheme with stepsize $T/n$ given by (\ref{eq:taylor-scheme}). In fact, doing so we obtain the numerical flow
\begin{equation*} \Psi(z;t_k,t_{k+1}):= z 
+   \sum_{i=1}^m \sigma^{(i)}(z)  \der B^{(i),n,T}_{t_k t_{k+1}} 
+   \sum_{i,j=1}^m  \mathcal{D}^{(i)} \sigma^{(j)}(z) \int_{t_k}^{t_{k+1}}\der B^{(i),n,T}_{t_k s}  \, d B_s^{(j),n,T} .\end{equation*}
 Since $B^{n,T}$ is the piecewise linear interpolation of $B$ on $[0,T]$ with stepsize $T/n$, the above iterated integrals can be now expressed as  products of increments of $B$. Indeed, according to the fact that
 \begin{equation}\label{eq:exp-B-n}
 \der B^{(i),n,T}_{t_k u}=\der B^{(i)}_{t_{k}t_{k+1}}  \frac{u-t_k}{T/n}, \quad
\dot{B}_u^{n,T} = \frac{n}{T} (\der B)_{t_{k}t_{k+1}} 
\quad\mbox{for}\quad u \in (t_k,t_{k+1}),
\end{equation}
it is readily checked that
\begin{equation*}
\der B^{(i),n,T}_{t_k t_{k+1}} = \der B^{(i)}_{t_{k}t_{k+1}},
\quad\mbox{and}\quad
\int_{t_k}^{t_{k+1}}\der B^{(i),n,T}_{t_k s}  \, d B_s^{(j),n,T}
=\frac{1}{2}\, \der B^{(i)}_{t_{k}t_{k+1}} \, \der B^{(j)}_{t_{k}t_{k+1}}.
\end{equation*}
Moreover, invoking relation (\ref{eq:exp-B-n}) and Lemma \ref{cont-fbm}, we get
$$   \left| \int_{t_k}^{t_{k+1}} |\dot{B}_u^{n,T}| \, du \right|   \leq 
\theta_{H,T} \, n^{-H} \, [\log(n)]^{1/2},
$$ 
for  $n$ large enough. Consequently, relation (\ref{ode-error-2-bis}) yields
\begin{align} \label{WZ-error-1b}
 \sup_{p,q=0, 1, \ldots, n-1, \, p \neq q} \frac{ | \der (Z^n- \overline{Z}^{n})_{t_p t_q} |    }{ |t_p-t_q|^{\ga}} \leq \theta_{H,\sigma, T} \, n^{-3H+1} \, [\log(n)]^{3/2}
\end{align}
for all $\gamma <H$ and all $n$ large enough.

This gives in particular
\begin{align} \label{WZ-error-1}
 \sup_{p,q=0, 1, \ldots, n-1, \, p \neq q} \frac{ | \der (Z^n- \overline{Z}^{n})_{t_p t_q} |    }{ |t_p-t_q|^{\ga}} \leq \theta_{H,\ga, \sigma, T} \, n^{-(H-\gamma)} \, [\log(n)]^{1/2}
\end{align}
for $1/3 < \gamma < H$.

\smallskip

Now it remains to "lift" this error estimate to $[0,T]$. For this we need the following smoothness result for the Wong-Zakai approximation.

 \begin{lemma}{\label{cont-wz}}
 Let $T>0$ and recall that $\overline{Z}^n$ is defined by equation (\ref{wz_2}). Then there exists $h^{**}>0$ and a finite and non-negative random variable $\theta_{H, h^{**},\sigma, T}$ such that for all $h\in (0,h^{**})$ and all $n\geq \frac{T}{h^{**}}$  we have
 $$ \sup_{t \in \lc 0, T- h \rc } |(\der \overline{Z}^n)_{t,t+h}| \leq \theta_{H, h^{**}, \sigma, T} \cdot h^H \cdot \sqrt{|\log(1/h)|} .$$
 
 \end{lemma}

\begin{proof}
As already mentioned in the proof of Lemma \ref{lemma_lin_int}, note that there exists $x_H>0$ such that the map $x \mapsto x^H \sqrt{|\log(1/x)|}$ is increasing on $(0,x_H]$. Set $h^{**}=\min(x_H,h^*)$, where $h^*$ is defined by Lemma \ref{cont-fbm}, and let $s,t\in [0,T]$ such that $\lln t-s \rrn \leq h^{**}$.

\smallskip

\noindent
(i) From (\ref{wz_2}) and (\ref{estim-remainder-term}), we deduce
\bean
\lln (\der \overline{Z}^n)_{st}- \sigma(\overline{Z}_s^n)(\der B^{n,T})_{st}\rrn &\leq &\lln t-s \rrn^{2\ka} G(\|  {\bf B}^{n,T} \|_{\ga})
\eean
for $1/3<\ka < \gamma < H$ and an increasing function $G: \R \rightarrow \R^+$.
Choosing $\ka,\gamma$ sufficiently large, we obtain
\bean \lln (\der \overline{Z}^n)_{st}- \sigma(\overline{Z}_s^n)(\der B^{n,T})_{st}\rrn 
&\leq & \theta_{H, h^*, \sigma, T} \lln t-s \rrn^H \sqrt{\log \lp \frac{1}{\lln t-s \rrn}\rp}.
\eean

 \noindent
(ii) Assume that $t_l \leq s \leq t \leq t_{l+1}$. One has
$$ \left| \sigma(\overline{Z}_s^n)(\der B^{n,T})_{st} \right| \leq \theta_{H, h^*, \sigma, T} \cdot  |t-s| \cdot (n/T)^{1-H} \sqrt{|\log(n/T)|}.$$
Since $ | t-s| \leq T/n $, i.e. $n/T \leq 1 / (t-s)$, it follows
$$ |(\der \overline{Z}^n)_{st}| \leq  \theta_{H, h^*, \sigma, T} \cdot (t-s)^H \cdot \sqrt{|\log(1/(t-s))|}.$$

\noindent
(iii) Now let $t_{l-1} \leq s \leq t_l  \leq t_p \leq t \leq t_{p+1}$ with $l \leq p$.
Then 
\begin{equation}\label{eq:dcp-iii}
(\der B^{n,T})_{st}  = (B_{t}^{n,T}-B_{t_{p}}) + (\der B)_{t_l t_{p}} + (B_{t_l} -B_{s}^{n,T}).
\end{equation}
As in the proof of Lemma \ref{lemma_lin_int}, this easily yields
\begin{equation}\label{eq:estim-iii}
\left| \sigma(\overline{Z}_s^n)(\der B^{n,T})_{st} \right| \leq   \theta_{H, h^{*}, \sigma, T}  \cdot (t-s)^H \cdot \sqrt{|\log(1/(t-s))|}
\end{equation}
for $|t-s| \leq T/n$. Whenever $|t-s| > T/n$, decomposition (\ref{eq:dcp-iii}) gives
$$ \left| \sigma(\overline{Z}_s^n)(\der B^{n,T})_{st} \right| \leq 2 \theta_{H, h^{*}, \sigma, T}   \cdot (T/n)^{H} \sqrt{|\log(n/T)|} + \theta_{H, h^{*}, \sigma, T} \cdot (t_p-t_l)^H \cdot \sqrt{|\log(1/(t_p-t_l))|} .$$
Using that $x \mapsto x^H \sqrt{|\log(1/x)|}$ is increasing, relation (\ref{eq:estim-iii}) is easily recovered.

\smallskip

\noindent
(iv) Combining the steps (i)-(iii) yields the assertion.

\end{proof}

\begin{proposition}\label{step_2} Let $T>0$ and $1/3 <  \gamma < H$. Then, there exists a finite and non-negative random variable $\eta_{H,\gamma, \sigma, T}^{(2)}$ such that
$$      \| Z^n- \overline{Z}^{n} \|_{\ga, \infty,T}   \leq  \eta_{H, \ga, \sigma, T}^{(2)} \cdot \sqrt{\log(n)} \cdot n^{-(H-\gamma)}$$ for $n>1$.
\end{proposition}

\begin{proof} Denote by $U^n$ the piecewise linear interpolation  with stepsize $T/n$ of the Wong-Zakai approximation $\overline{Z}^n$. Proceeding as in the proof of 
Lemma \ref{lemma_lin_int} and using Lemma \ref{cont-wz} we have
$$     \| U^n- \overline{Z}^{n} \|_{\ga, \infty,T}   \leq  \theta_{H, \gamma, \sigma, T} \cdot \sqrt{\log(n)} \cdot n^{-(H-\gamma)}.  $$
Thus, it remains to consider the difference between $U^n$ and $Z^n$. For $t \in \lc t_k , t_{k+1} \rc$ for some $k$ we have
$$ U^n_t- Z^n_t=  \overline{Z}^{n}_{t_k}-{Z}^{n}_{t_k} + \frac{t-t_k}{T/n} \der\lp \overline{Z}^{n}-{Z}^{n}\rp_{t_kt_{k+1}}.$$
Assuming additionally that $s \in \lc t_{l}, t_{l+1} \rc$ and $t\in [t_k,t_{k+1}]$ for some $l \leq k$, we have
\begin{equation}\label{eq:der-U-Z}
\der (U^n- Z^n)_{st} =  \der(\overline{Z}^{n}-{Z}^{n})_{t_lt_k} 
 + \frac{t-t_k}{T/n} \der (\overline{Z}^{n}-{Z}^{n})_{t_k t_{k+1}}  
  - \frac{s-t_l}{T/n} \der (\overline{Z}^{n}-{Z}^{n})_{t_l t_{l+1}}   .
\end{equation}

\noindent
(i) Assume that $l+1 <k$.
 Applying (\ref{WZ-error-1}) to relation (\ref{eq:der-U-Z}) and according to the fact that $(s-t_{l}) \leq T/n$, $(t-t_{k})\le T/n$, we obtain 
 \begin{equation}\label{eq:estim-der-U-Z}
|\der (U^n- Z^n)_{st}| \leq  
\theta_{ H,  \ga, \sigma, T}  \,  |t-s|^{\ga} \cdot n^{-(H-\ga)} \sqrt{\log(n)} .
\end{equation}

 \noindent
 (ii) Assume that $l=k$.
 Here (\ref{eq:der-U-Z}) simplifies to
 \begin{align*}
 \der (U^n- Z^n)_{st}&=   \frac{t-s}{T/n} \der (\overline{Z}^{n}-{Z}^{n})_{t_kt_{k+1}}  \end{align*}
and thus (\ref{WZ-error-1}) combined with the fact that $|t-s|\le T/n$ gives an estimate of the form~(\ref{eq:estim-der-U-Z}) again.

\smallskip

Finally, the case $k=l+1$ can be treated in a similar manner, and this completes the proof.

\end{proof}

\begin{remark}
Putting together Propositions \ref{step_1} and \ref{step_2}, our Main Theorem \ref{main_thm} now follows.
\end{remark}

\smallskip

\subsection{Optimality of the error bound}\label{optimality}
Reviewing the steps of the derivation of our main result, one realises that the final convergence rate $n^{-(H-\ga)} \sqrt{\log (n)}$ is directly linked to the error (measured in the $\ga$-H\"older norm) of the piecewise linear interpolation of fractional Brownian motion. All other estimates lead to higher rates of convergence. As a result, in order to prove the optimality of our result, it is natural to consider the most simple equation
$$ dY^{(1)}_t=dB^{(1)}_t, \quad t\in \lc 0,T \rc, \qquad Y_0=a\in \R,$$
for which our Milstein-type approximation  is given by $Z^n=B^{n,T}$.

\smallskip

First, observe that
\begin{align*}
   \| Y- Z^{n} \|_{\gamma, \infty,T}    =  \| B^{(1)}- B^{(1),n,T} \|_{\gamma, \infty,T}  &  \geq  \sup_{s,t \in [0,T]} \frac{ | \der(B^{(1)} -B^{(1),n,T})_{st}| } {|t-s|^{\gamma} }
    \\ & \geq  \sup_{t \in [0,T]} \frac{ | B_t^{(1)} - B_{t}^{(1),n,T}| } {t^{\gamma} }.
\end{align*}
Using the scaling and stationarity properties of fBm, we get
\begin{multline}\label{optimality-proof}
   \sup_{t \in [0,T]} \frac{ | B_t - B_{t}^{n,T}| } {t^{\gamma} } \stackrel{\mathcal{L}}{=}  
   \sup_{t \in [0,1]}  T^H\frac{ | B_t - B_{t}^{n,1}| } {T^{\ga} t^{\gamma} } 
   \stackrel{\mathcal{L}}{=} T^{ H - \ga}
   \sup_{t \in [0,n]} n^{-H}  \frac{ | B_t - B_{t}^{n,n}| } { n^{-\ga} t^{\gamma} } \\  
    \geq n^{-(H-\gamma)}  T^{ H-\ga}\sup_{t \in [1,n]}  | B_t - B_{t}^{n,n}| \stackrel{\mathcal{L}}{=}   n^{-(H-\gamma)} T^{H- \ga}\sup_{t \in [0,n-1]} | B_t - B_{t}^{n-1,n-1}|   
  . \end{multline}
  
Now let us recall the following result of \cite{huesleretal}:
 $$\frac{v_n}{\si_n} \lp \sup_{t\in [0,1]} \lln B_t-B^{n,1}_t \rrn -\si_n v_n \rp \stackrel{\mathcal{L}}{\longrightarrow}  G,$$
    where $G$ is a Gumbel distribution, $\lim_{n \rightarrow \infty} \frac{v_n}{\sqrt{2  \log(n)}} =1$ and $ \lim_{n \rightarrow \infty} n^{H}\si_n =c_H$. This implies in particular
$$  \frac{n^H}{\sqrt{ \log(n)} }  \sup_{t\in [0,1]} \lln B_t-B^{n,1}_t \rrn       \stackrel{Prob.}{\longrightarrow}   \sqrt{2} c_H.$$
 Applying again the scaling property of fBm gives  
$$ \frac{1}{\sqrt{ \log(n)}}   \sup_{t\in [0,n]} \lln B_t-B^{n,n}_t \rrn   \stackrel{\mathcal{L}}{\longrightarrow} \sqrt{2} c_H $$  and so
   $$ \frac{1}{\sqrt{ \log(n)}}   \sup_{t\in [0,n-1]} \lln B_t-B^{n-1,n-1}_t \rrn   \stackrel{\mathcal{L}}{\longrightarrow} \sqrt{2} c_H.$$  

Going back to (\ref{optimality-proof}), this finally yields

$$       \lim_{n \rightarrow \infty}      {\bf P} \left(\,  \ell(n) \cdot \| Y- Z^{n} \|_{\gamma, \infty,T} \, < \, \infty \right)=0  , $$
 if $$\liminf_{n \rightarrow \infty}  \, \ell(n) \cdot \sqrt{\log(n)} \cdot n^{-(H-\gamma)} = \infty,$$
which corresponds to our claim at Remark \ref{rem_low_bound}.

\medskip

\section{Numerical Examples}

In the introduction, we stated the conjecture that the error in the supremum norm of our proposed modified Milstein scheme satisfies
$$ \| Y - Z^{n} \|_{\infty, T}  \leq \eta_{ H, \sigma, T} \cdot  \sqrt{\log(n)} \cdot \big( n^{-H} + n^{-2H+1/2} \big).$$ 
\begin{figure}[htp]
\epsfig{figure=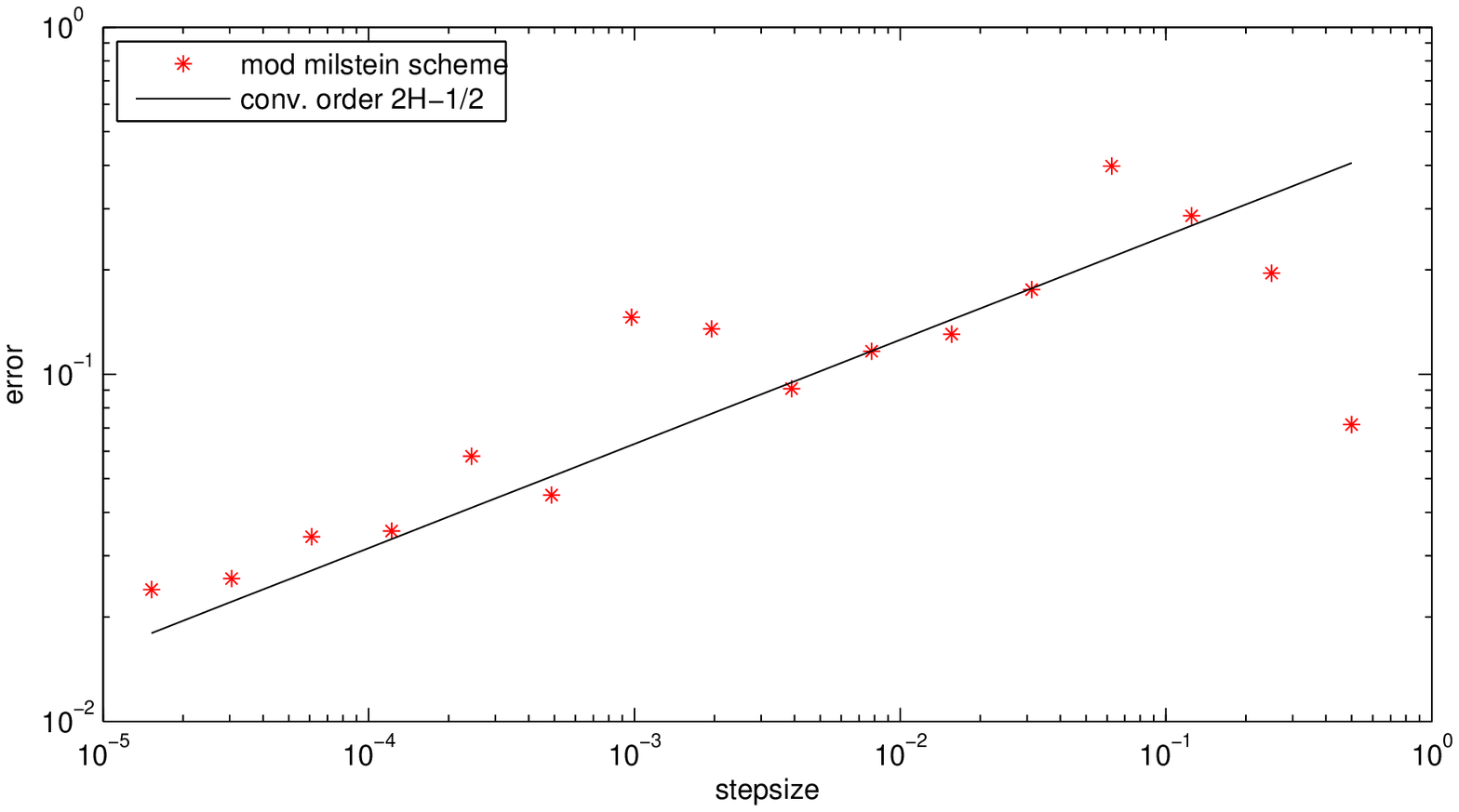 ,width=6.5cm, height=5.5cm}
\ \leavevmode
\epsfig{figure=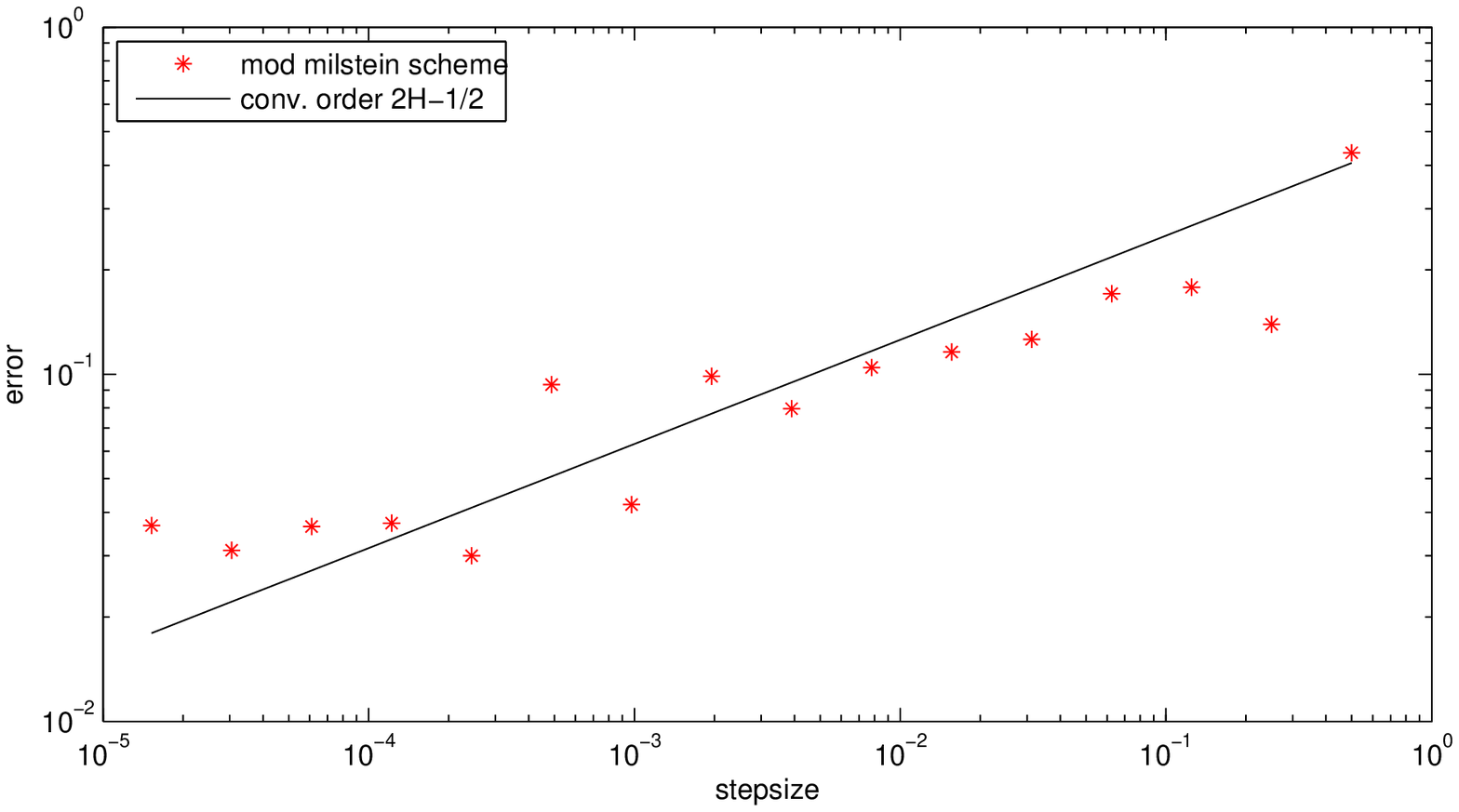, width=6.5cm, height=5.5cm}
\ \leavevmode
\epsfig{figure=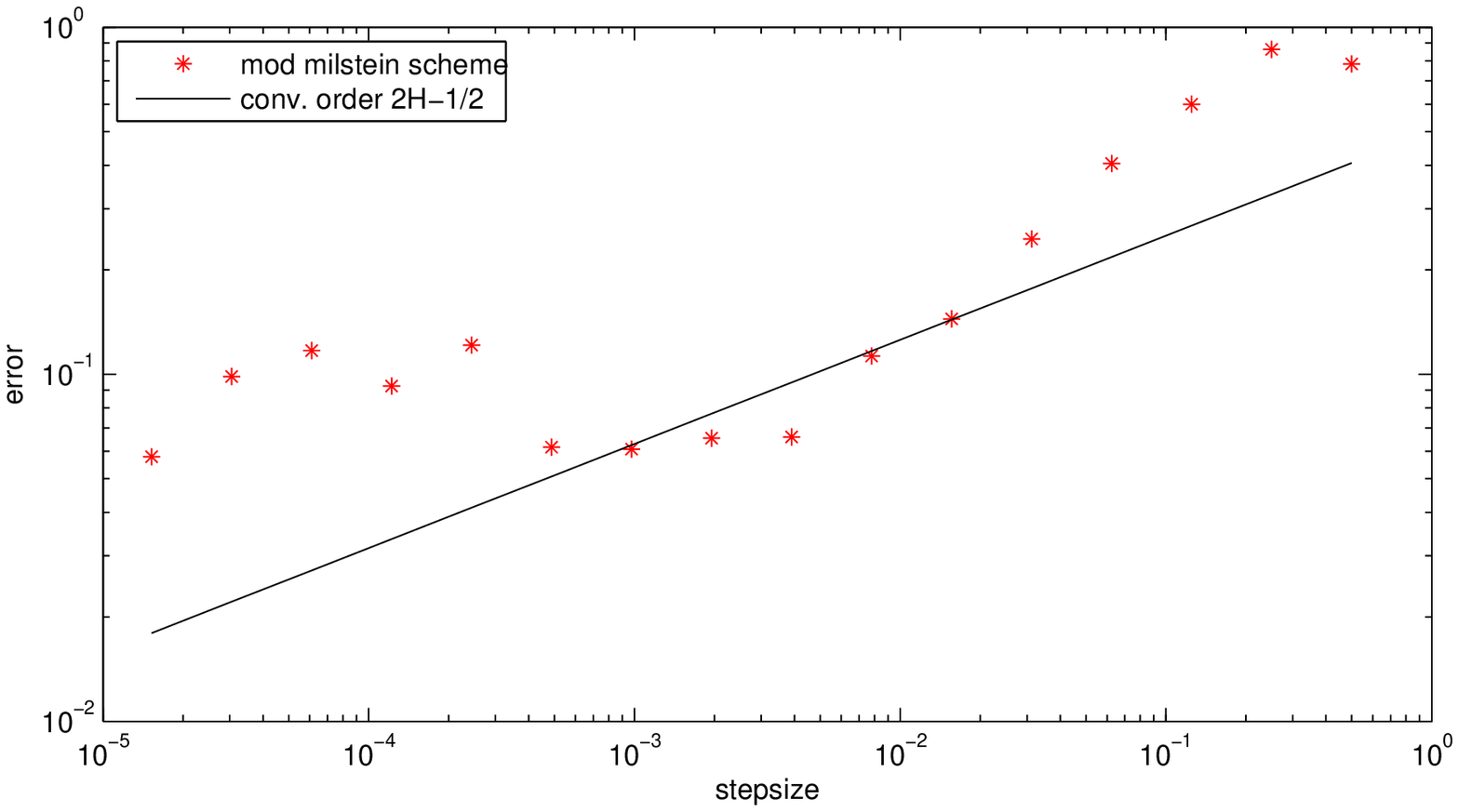, width=6.5cm, height=5.5cm}
\ \leavevmode
\epsfig{figure=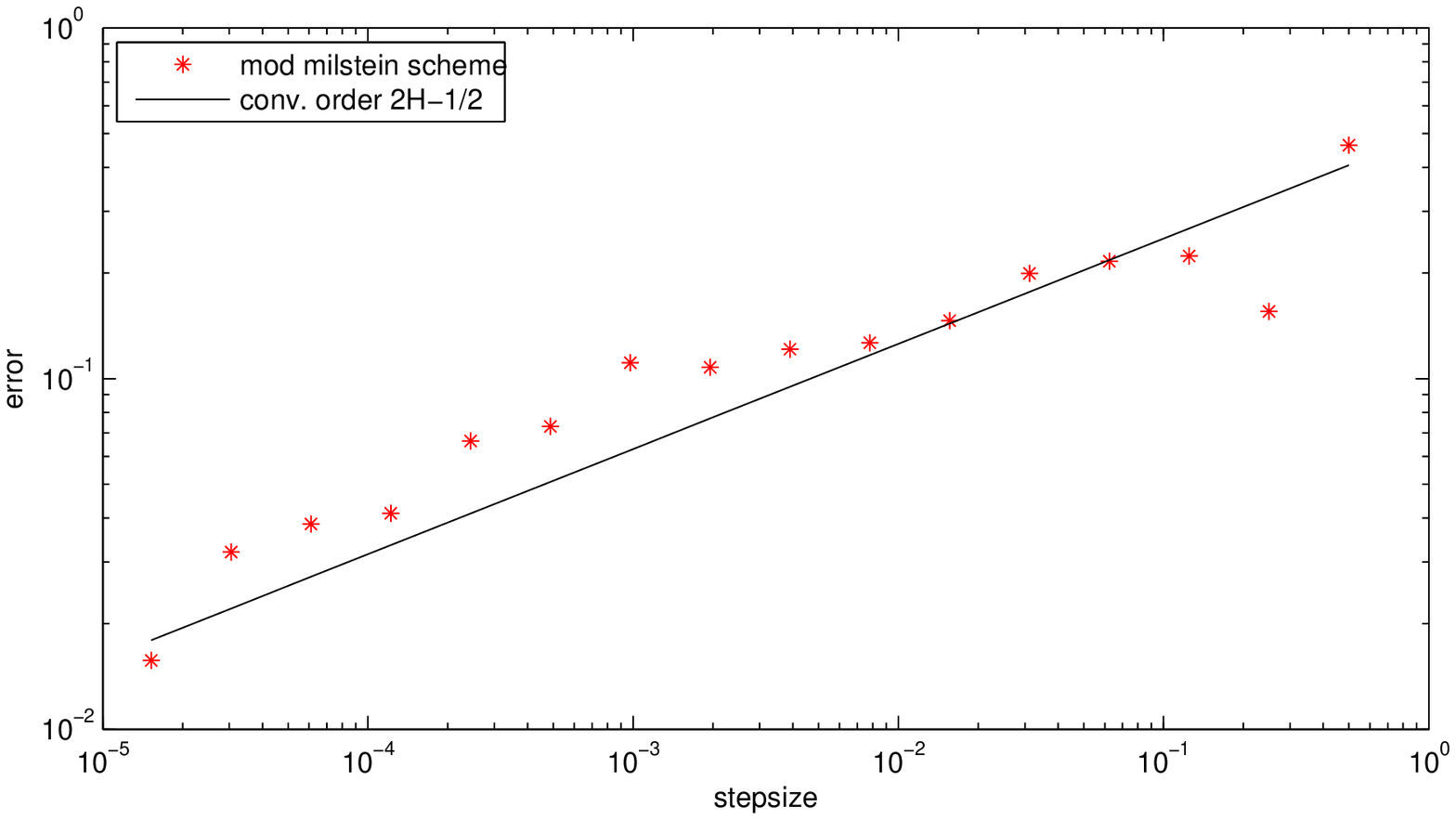, width=6.5cm, height=5.5cm}
   \caption{Equation (\ref{test_bound_coeff}):  pathwise maximum error  vs. step size for four
   sample paths for   $H=0.4$. }
\end{figure}
\begin{figure}[htp]
\epsfig{figure=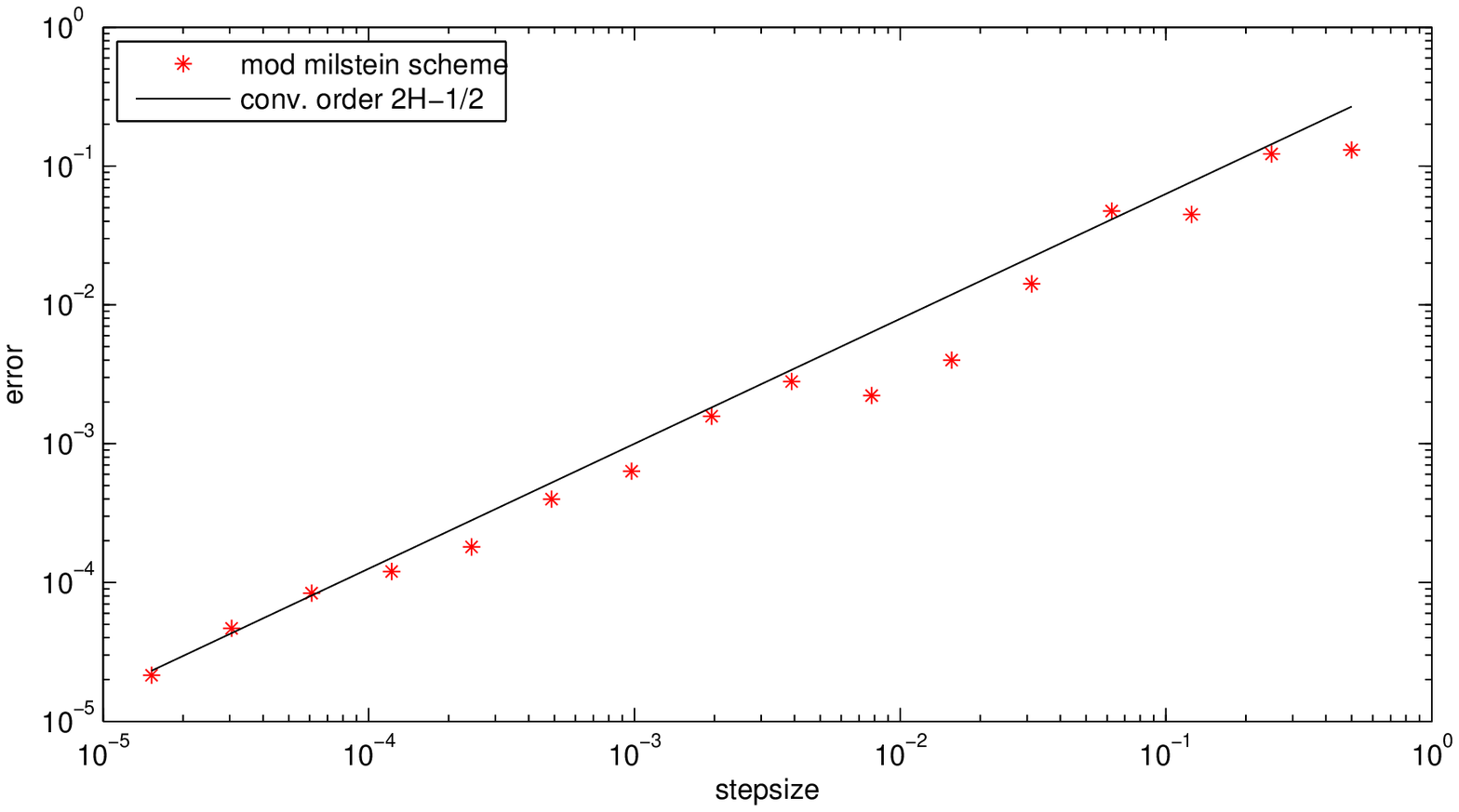 ,width=6.5cm, height=5.5cm}
\ \leavevmode
\epsfig{figure=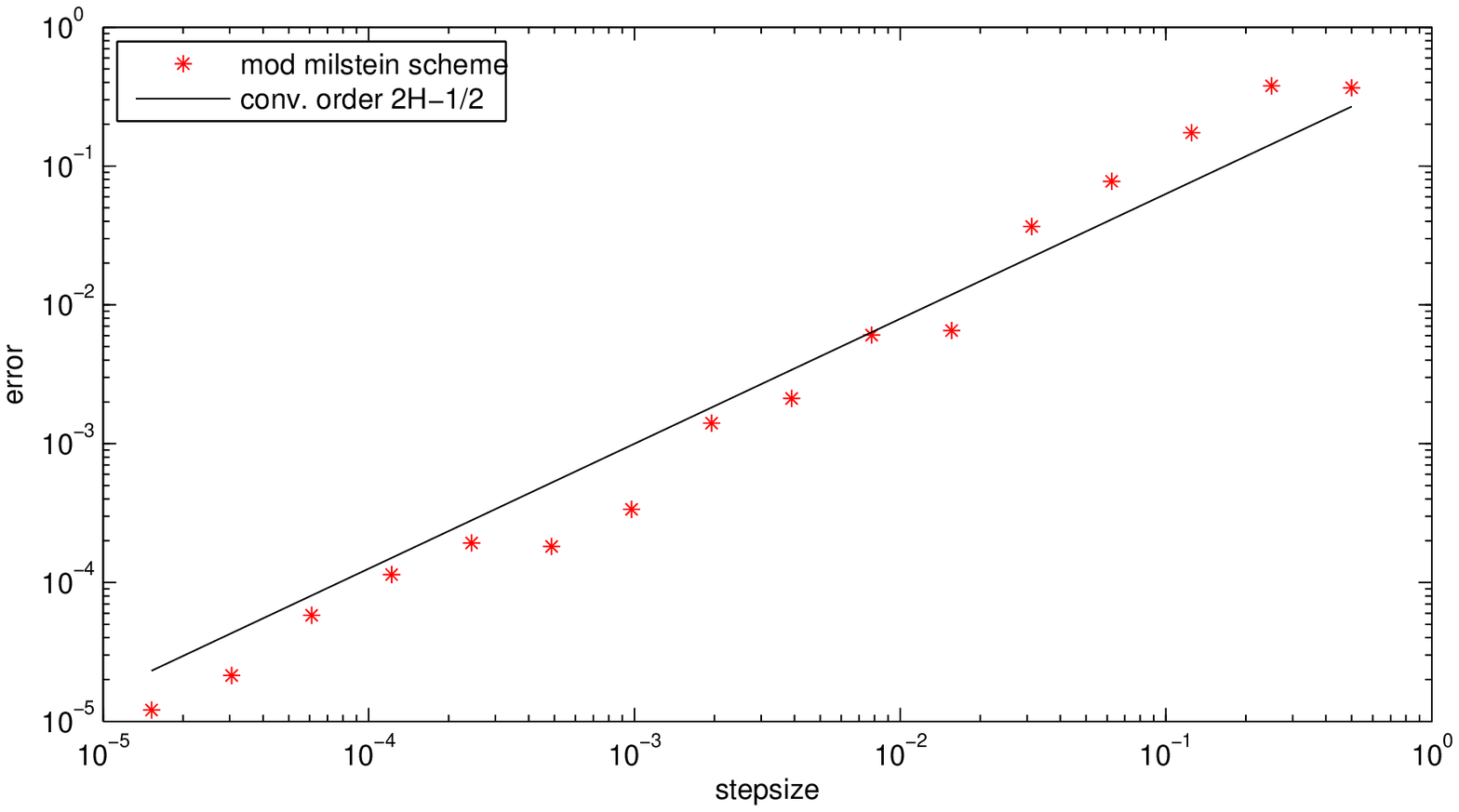, width=6.5cm, height=5.5cm}
\ \leavevmode
\epsfig{figure=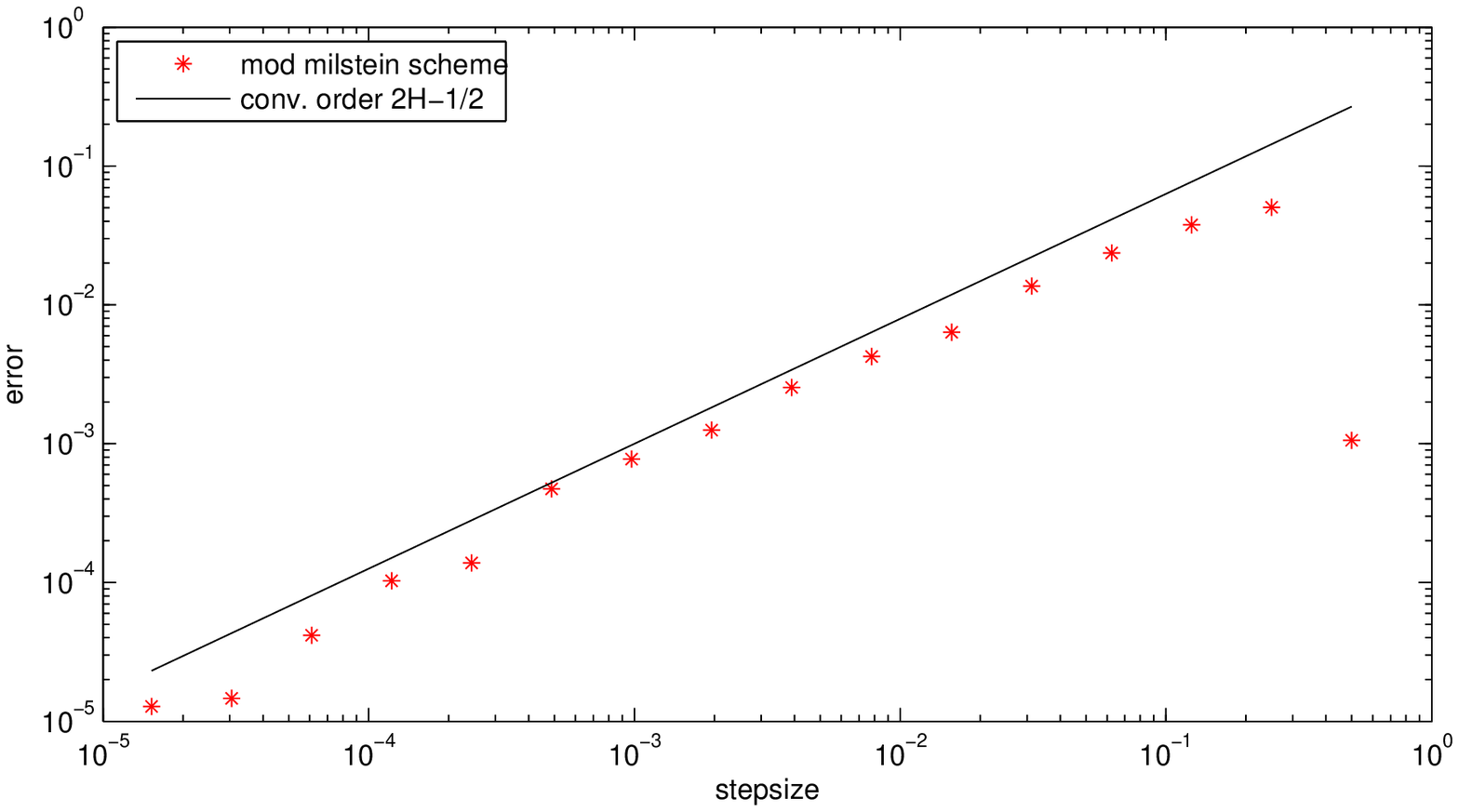, width=6.5cm, height=5.5cm}
\ \leavevmode
\epsfig{figure=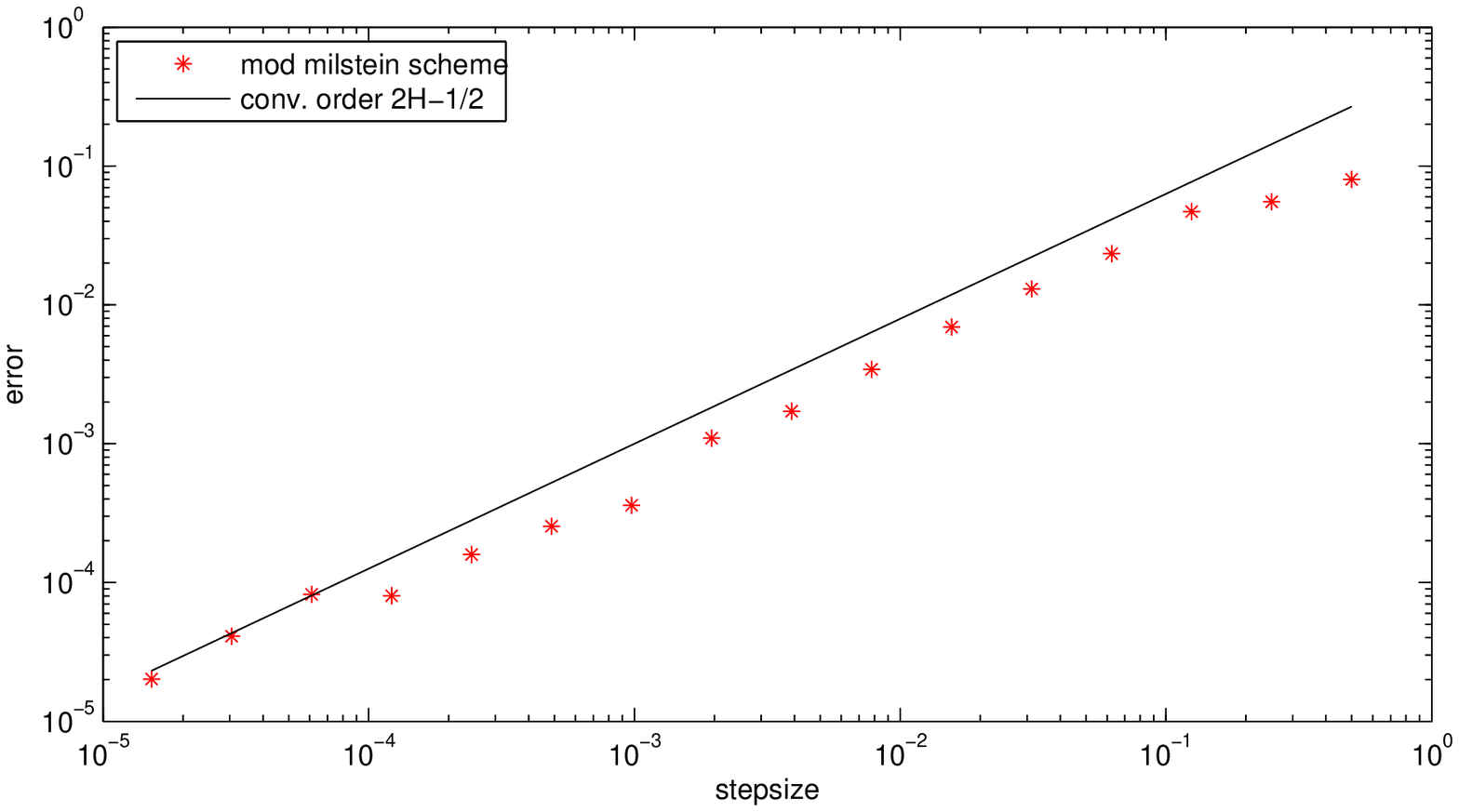, width=6.5cm, height=5.5cm}
   \caption{Equation (\ref{test_bound_coeff}):  pathwise maximum error  vs. step size for four
   sample paths for    $H=0.7$. }
\end{figure}

Note that if $U^n$ denotes the piecewise linear interpolation of $Z$ with stepsize $T/n$, then we have
$$ \| Y - U^{n} \|_{\infty, T}  \leq \eta_{ H, \sigma, T}  \cdot  \sqrt{\log(n)} \cdot  n^{-H},$$ 
which follows from a straightforward modification of the Lemmata \ref{lemma_lin_int} and \ref{cont-wz}.
Since furthermore
$$ \| Y - Z^{n} \|_{\infty, T}  \leq  \| Y - U^{n} \|_{\infty, T}  + \max_{k=0, \ldots, n} | Y_{kT/n} - Z^{n}_{kT/n} |, $$
it suffices to consider the maximal error in the discretisation points, i.e.
$$ \max_{k=0, \ldots, n} | Y_{kT/n} - Z^{n}_{kT/n} |, $$ 
to support our conjecture.

Our first example will be the SDE
\begin{align} dY_t = \cos(Y_t) \, dB^{(1)}_t + \sin(Y_t)\,d B^{(2)}_t, \quad t \in [0,1] , \qquad Y_0=1. \label{test_bound_coeff} \end{align}
Figure 1 shows the maximum error in the discretization points,  
i.e.  $$ \max_{k=0, \ldots, n} | Y_{kT/n} (\omega)- Z^{n}_{kT/n}(\omega) |, $$ which for brevity we call in the following maximum error,  versus the
step size  $1/n$ for four  different sample paths $\omega  \in \Omega$ for $H=0.4$,  while  Figure 2 shows the maximum error versus the
step size $1/n$ for four  different sample paths $\omega  \in \Omega$ for $H=0.7$. (So small values on the $x$-axis correspond to small stepsizes, while small values on the $y$-axis correspond to small errors and vice versa.)

The numerical reference solution is obtained by using our Milstein-type scheme with very small stepsize.
Since we use 
log-log-coordinates, the straight lines correspond to the
convergence order $2H-1/2$.  The  stars  correspond to the error of the Milstein-type scheme.
For $H=0.4$  the estimated  convergence rates are  in acceptable accordance
with our conjecture, while for $H=0.7$ they are in good accordance.

\begin{figure}[htp]
\epsfig{figure=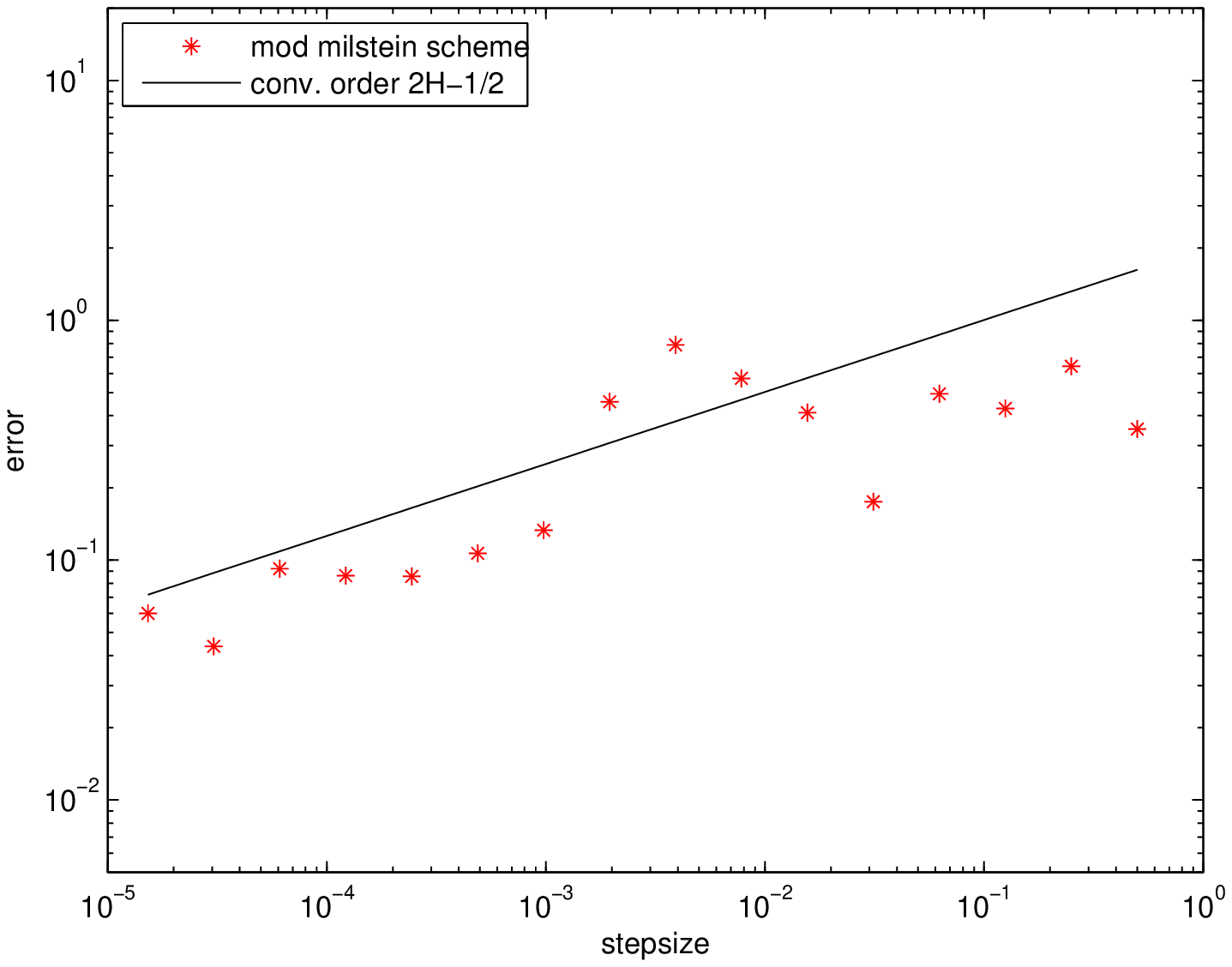 ,width=6.5cm, height=5.5cm}
\ \leavevmode
\epsfig{figure=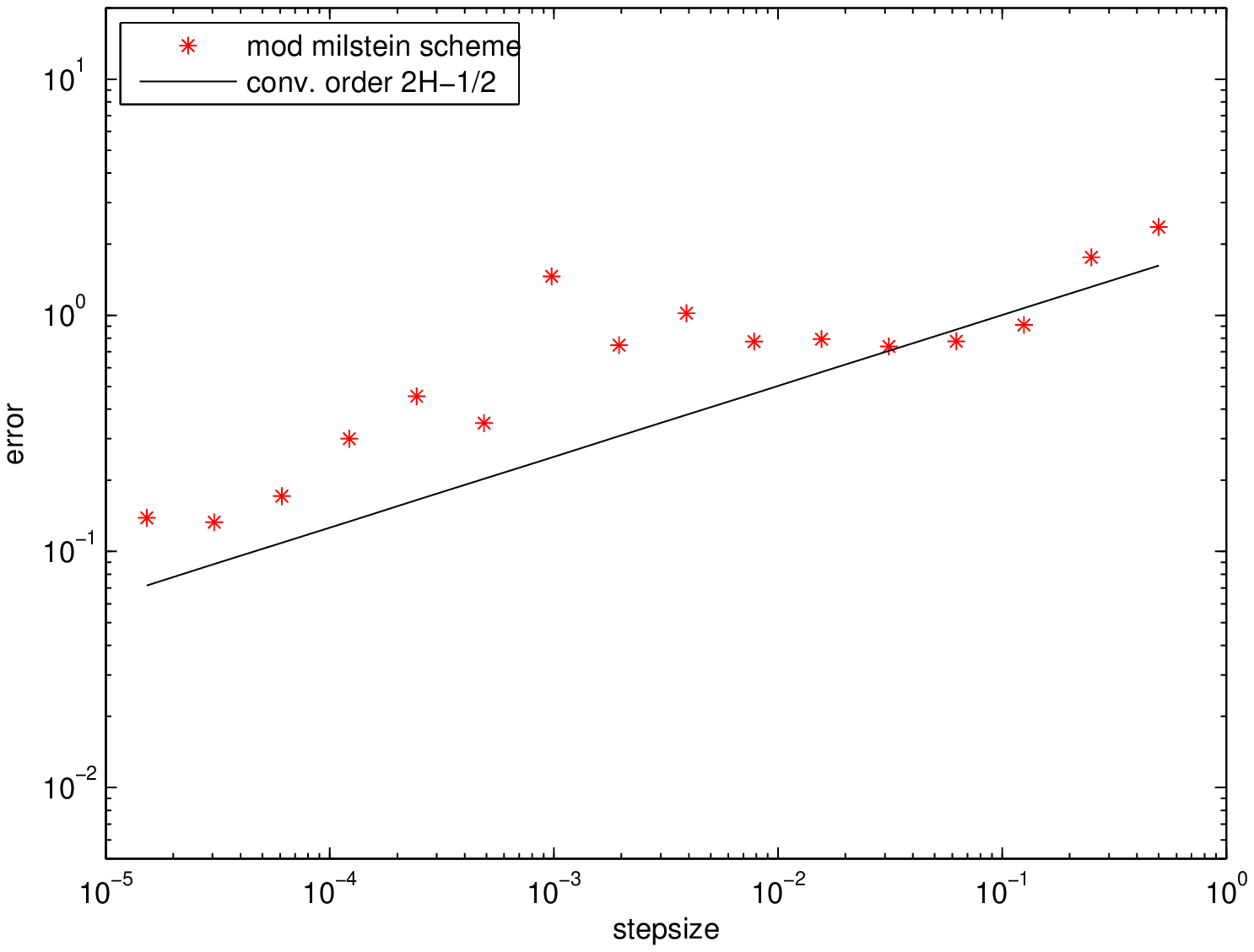, width=6.5cm, height=5.5cm}
\ \leavevmode
\epsfig{figure=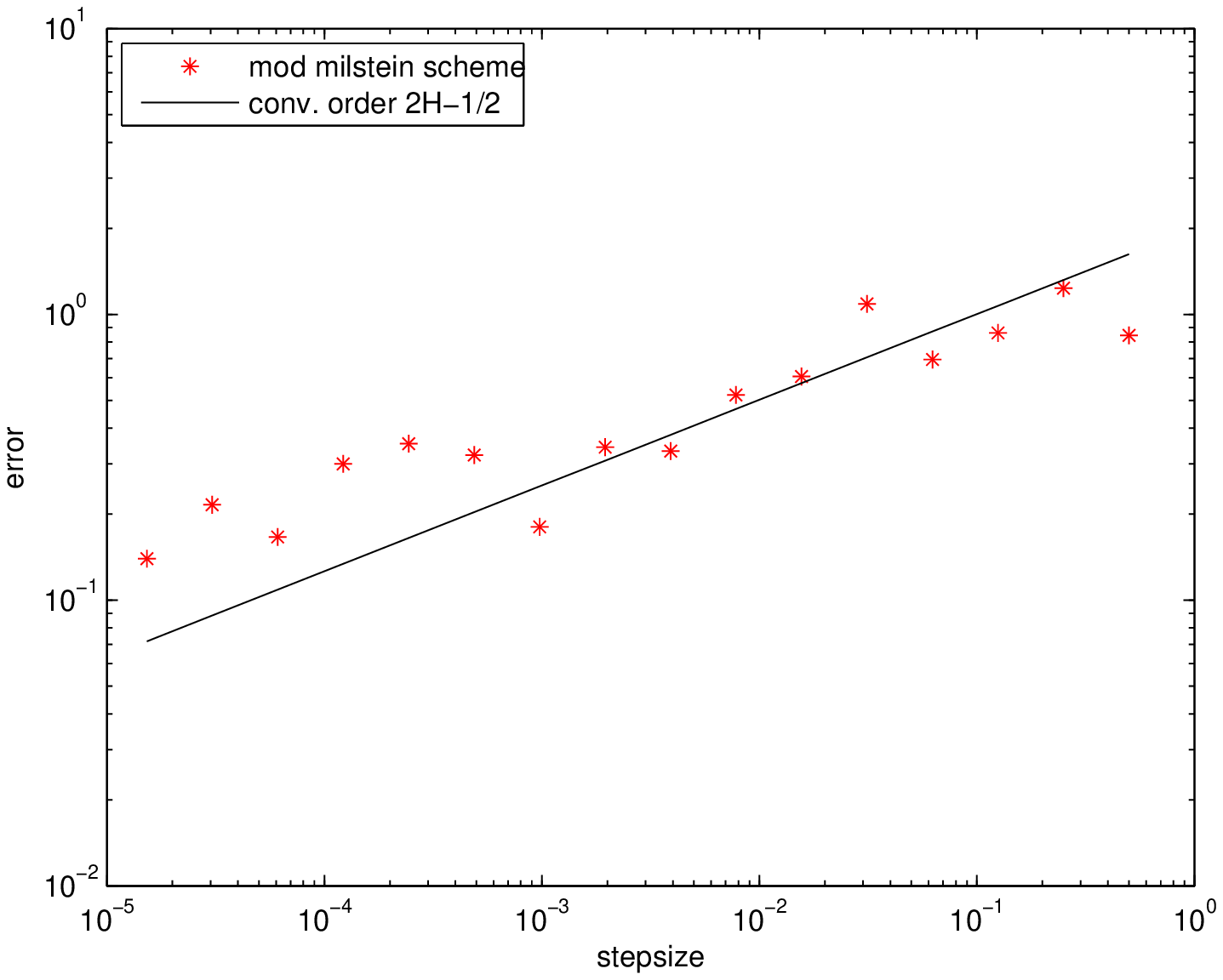, width=6.5cm, height=5.5cm}
\ \leavevmode
\epsfig{figure=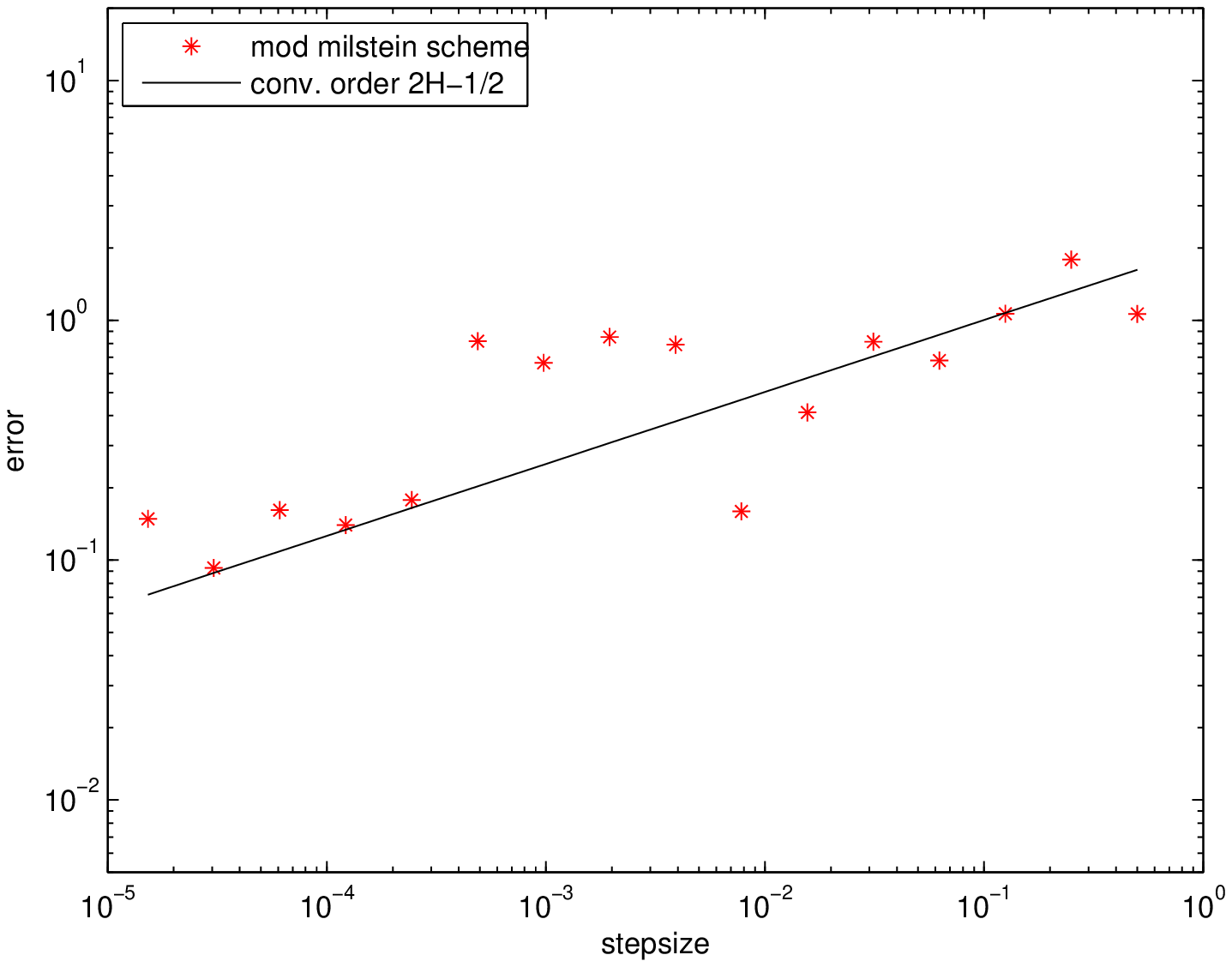, width=6.5cm, height=5.5cm}
   \caption{Equation (\ref{lin_coeff_test}):  pathwise maximum error  vs. step size for four
   sample paths for    $H=0.4$. }
\end{figure}
\begin{figure}[htp]
\epsfig{figure=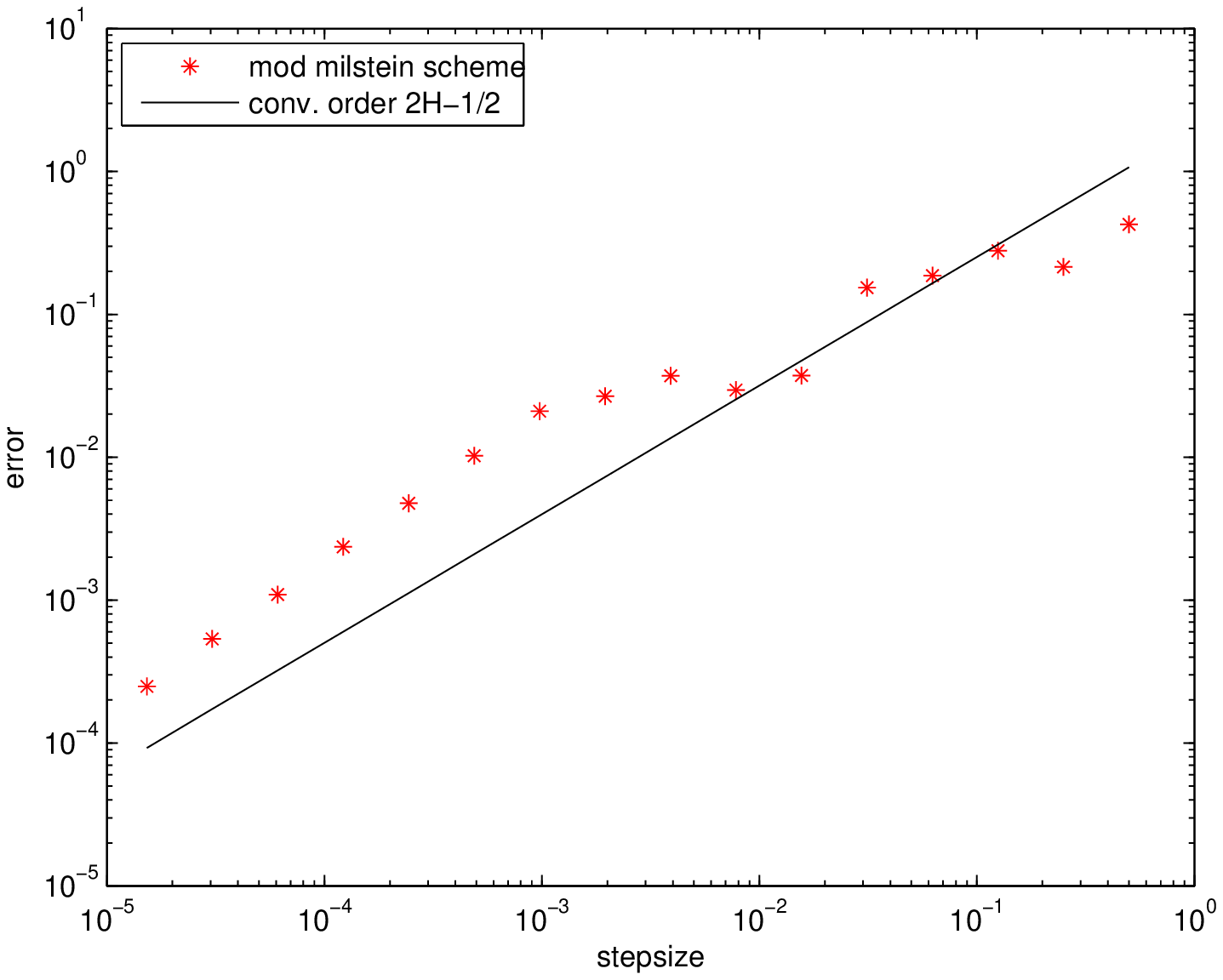 ,width=6.5cm, height=5.5cm}
\ \leavevmode
\epsfig{figure=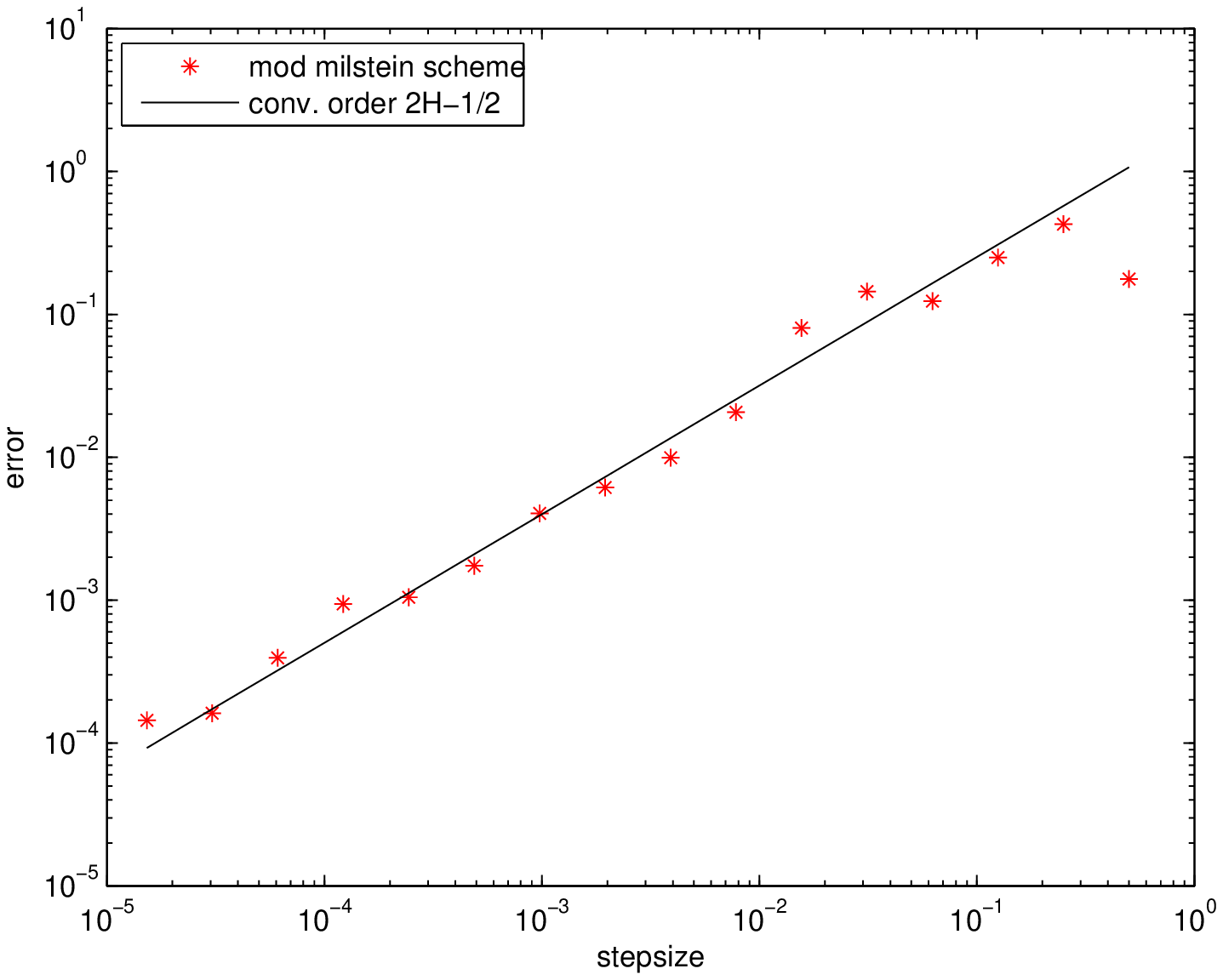, width=6.5cm, height=5.5cm}
\ \leavevmode
\epsfig{figure=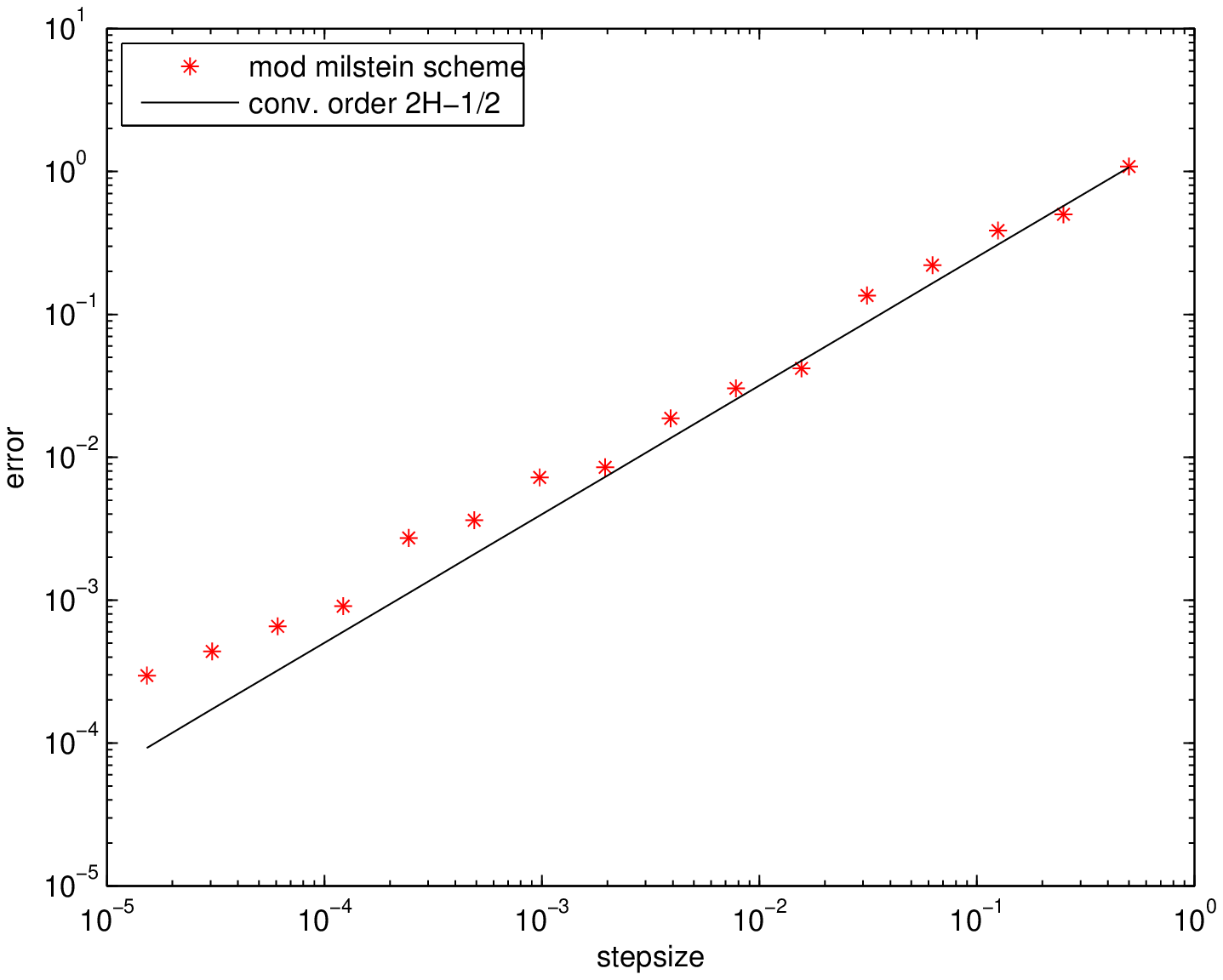, width=6.5cm, height=5.5cm}
\ \leavevmode
\epsfig{figure=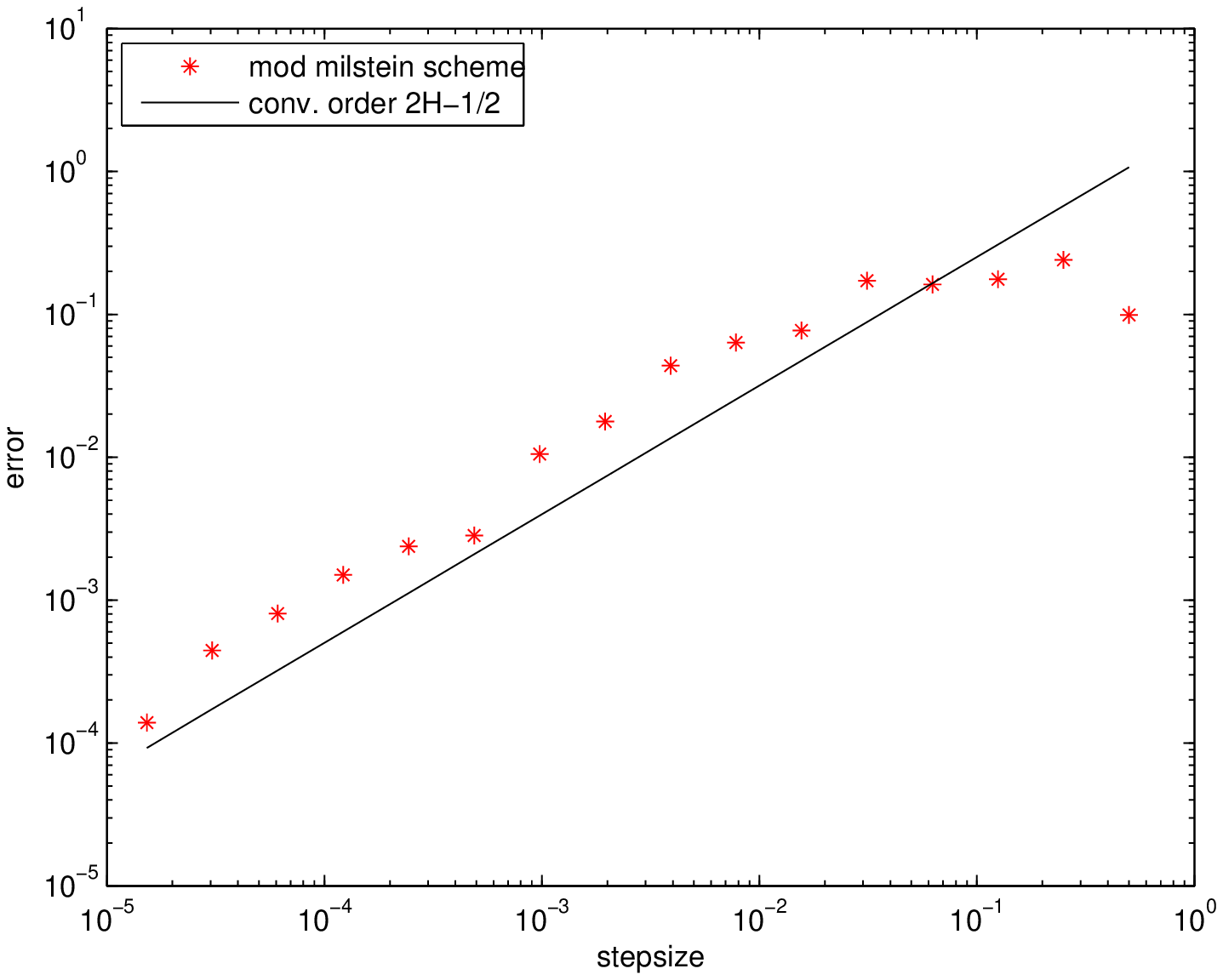, width=6.5cm, height=5.5cm}
   \caption{Equation (\ref{lin_coeff_test}):  pathwise maximum error  vs. step size for four
   sample paths for    $H=0.7$. }
\end{figure}

As second example we consider the linear equation
\begin{align}
dY_t^{(1)}= Y_t^{(2)}dB_t^{(1)}, \quad dY_t^{(2)}= Y_t^{(1)}dB_t^{(2)} ,\quad t \in [0,1], \qquad \,\, Y_0^{(1)}= 1, \, Y_0^{(2)}=2 \label{lin_coeff_test}. \end{align}
Figures 3 and 4  show again  the maximum error versus the step size   for four  different sample paths for $H=0.4$ and $H=0.7$, respectively. Again the estimated covergence rates are in acceptable accordance  with our conjecture for $H=0.4$ and in  good accordance for $H=0.7$. 

Note that the convergence order $2H-1/2$ is quite slow for  small $H$. In particular, for $H=0.4$ the convergence order equals $0.3$. We suppose that this effect also causes the fluctuating behaviour in the
estimated convergence rates in the case $H=0.4$.

\bigskip

\section{Appendix: Proof of Theorem \ref{thm:Lip}}

\subsection{Existence and uniqueness of the solution}
This section gives some details of the proof of point (1)  of Theorem 
\ref{thm:Lip} in the case $\gamma \leq 1/2$. The case $\gamma >1/2$ is simpler and thus omitted. 

\smallskip

The solution to equation (\ref{de}) is obtained via a fixed-point argument, which is first applied locally and then extended to the whole interval $[0,T]$.
\smallskip

\noindent
\textit{Notations}. For $\cq_{\ka,a}^x([\ell_1, \ell_2];\R^d)$ we will write in the following only $\cq_{\ka}^x([\ell_1,\ell_2])$ to simplify the notation. In particular, note that the norm $\cn[ \cdot;\cq_{\ka,a}^x([\ell_1,\ell_2])]$ does not depend on $a \in \R^d$. Moreover, for $y \in \cq_{\ka}^x([\ell_1,\ell_2])$, which admits the decomposition
$$ (\delta y)_{st} = \zeta_s  (\delta \xi)_{st} + r_{st},$$
we  set 
$$ y^x:= \zeta, \qquad y^{\sharp}:=r. $$

\smallskip

\noindent
\textit{Local considerations}. Consider a time $0<T_0 \leq T$ and for any $y\in \cq_\ka^x([0,T_0])$, define $z=\Gamma_{T_0}(y)$ as the unique process in $\cq_\ka^x([0,T_0])$ such that $z_0=y_0$ and $(\der z)_{st}=\cj_{st}(\si(y) \, dx)$. If $y,\yti \in \cq_\ka^x([0,T_0])$ with $(y_0,y^x_0)=(\yti_0,\yti^x_0)=(a,\si(a))$, and if $z=\Gamma_{T_0}(y),\zti=\Gamma_{T_0}(\yti)$, then some standard differential calculus easily leads to
\begin{equation}\label{relation-un}
\cn[z;\cq_\ka^x([0,T_0])] \leq c_x \lcl 1+T_0^{\ga-\ka} \cn[y;\cq_\ka^x([0,T_0])]^2 \rcl,
\end{equation}
and
\begin{multline}\label{relation-deux}
\cn[z-\zti;\cq_{\ka}^x([0,T_0])] \\
\leq c_x T_0^\ka \cn[y-\yti;\cq_{\ka}^x([0,T_0])] \lcl 1+\cn[y;\cq_\ka^x([0,T_0])]^2+\cn[\yti;\cq_\ka^x([0,T_0])]^2 \rcl,
\end{multline}
with $c_x=c (1+\norm{x}_\ga+\norm{\mathbf{x}^2}_{2\ga})$ for some constant $c>1$. Now set $T_0=(4c_x^2)^{-1/(\ga-\ka)}$ and $R_{T_0}=2c_x$, so that, if in addition $\cn[y;\cq_\ka^x([0,T_0])]\leq R_{T_0}$, then by (\ref{relation-un}), $\cn[z;\cq_\ka^x([0,T_0])]$ $\leq R_{T_0}$ and, if also  $\cn[\tilde{y};\cq_\ka^x([0,T_0])]\leq R_{T_0}$, by (\ref{relation-deux}),
$$\cn[z-\zti;\cq_{\ka}^x([0,T_0])] \leq c_x \, \cn[y-\yti;\cq_{\ka}^x([0,T_0])] \cdot (4c_x^2)^{-\ka/(\ga-\ka)}\lcl 1+8c_x^2 \rcl.$$
Observe that $3-2\ka/(\ga-\ka) <0$ for $1/3 < \kappa < \gamma \leq 1/2$ and so
$$c_x (4c_x^2)^{-\ka/(\ga-\ka)}\lcl 1+8c_x^2 \rcl =\lp \frac{1}{4}\rp^{\ka/(\ga-\ka)} \lcl c_x^{1-2\ka/(\ga-\ka)}+8c_x^{3-2\ka/(\ga-\ka)} \rcl \leq 9 \lp \frac{1}{4}\rp^2 <1.$$
As a result, $\Gamma_{T_0}$ is a strict contraction of the following closed subset of $\cq_{\ka}^x([0,T_0])$:
$$\cb_{(a,\si(a)),R_{T_0}}^{T_0}= \lcl y \in \cq_{\ka}^x([0,T_0]); \ (y_0,y^x_0)=(a,\si(a)) \ , \ \cn[y;\cq_\ka^x([0,T_0])] \leq R_{T_0} \rcl.$$

Let us denote by $y^{T_0}$ the fixed point of the restriction of $\Gamma_{T_0}$ to $\cb_{(a,\si(a)),R_{T_0}}^{T_0}$. 

\smallskip

\noindent
\textit{Extending the solution.} One can use the same arguments as in the previous step for the set 
\begin{multline}
\cb_{\lp y_{T_0}^{T_0},\si\lp y_{T_0}^{T_0}\rp \rp,R_{T_0}}^{2T_0}\\
= \lcl y \in \cq_{\ka}^x([T_0,2T_0]); \ (y_{T_0},y^x_{T_0})=\lp y_{T_0}^{T_0},\si\lp y_{T_0}^{T_0}\rp \rp \ , \ \cn[y;\cq_\ka^x([T_0,2T_0])] \leq R_{T_0} \rcl,
\end{multline}
and this provides us with an extension of the solution on $[T_0,2T_0]$, denoted by $y^{2T_0}$. Repeat the procedure until $[0,T]$ is covered, and then define
$$y=\sum_{i=1}^{N_{T_0}} y^{iT_0} \cdot 1_{[(i-1)T_0,iT_0]} \ , \ y^x=\sum_{i=1}^{N_{T_0}} y^{x,iT_0} \cdot 1_{[(i-1)T_0,iT_0]},$$
where $N_{T_0}$ is the smallest integer such that $N_{T_0} \cdot T_0 \geq T$.

\smallskip

It is not hard to see that $y$ is a solution to the system (\ref{de}). Moreover,
\begin{align*}
&\cn[y;\cq_\ka^x([0,T])]\\
&\leq  \sup_{k=1,\ldots, N_{T_0}} \cn[y^{kT_0};\cq_\ka^x([(k-1)T_0,kT_0])]
+\lcl 1+\norm{x}_\ga \rcl \sum_{k=1}^{N_{T_0}} \cn[y^{kT_0};\cq_\ka^x([(k-1)T_0,kT_0])]\\
&\leq  R_{T_0}+R_{T_0} \cdot N_{T_0} \cdot \lcl 1+\norm{x}_\ga \rcl
\leq  2c_x \lp 1+\lp T/T_0+1 \rp \lp 1+\norm{x}_\ga \rp \rp\\
&\leq  2c_x \lp 1+\lp 1+4\cdot T\cdot c_x^{2/(\ga-\ka)}\rp \lp 1+\norm{x}_\ga \rp \rp,
\end{align*}
which gives the estimate (\ref{control-sol-y}). The unicity of this solution is easy to prove due to (\ref{relation-deux}). The details are left to the reader.

\subsection{Continuity of the It\^o map}
We shall now prove point (2) in Theorem \ref{thm:Lip}. For this, let us again introduce some notation:

\smallskip

\noindent
\textit{Notation}: If $y\in \cq_\ka^x$ and $\yti \in \cq_\ka^{\xti}$ for two different driving signals $x,\xti$, define
\begin{multline*}
\cn[y-\yti;\cq_\ka^{x,\xti}]=\cn[(y,y^x)-(\yti,\yti^x);\cq_\ka^{x,\xti}]:=\cn[y-\yti;\cac_1^\ga]+\cn[y^x-\yti^x;\cac_1^{0,\ka}]\\
+\cn[y^\sharp-\yti^\sharp;\cac_2^{2\ka}].
\end{multline*}

\smallskip

\noindent
\textit{Local considerations}. Consider a time $T_0>0$. From the decomposition
\begin{multline*}
\der (y-\yti)_{st}=\lc \si(y_s)-\si(\yti_s)\rc \cdot (\der x)_{st}+\si(\yti_s) \cdot \der (x-\xti)_{st}+\lc y^x_s \si'(y_s)-\yti^x \si'(\yti_s)\rc \cdot \mathbf{x}^2_{st}\\ {}\qquad \qquad
+\yti^x_s \si'(\yti_s) \cdot \lc \mathbf{x}^2_{st}-\tilde{\mathbf{x}}^2_{st} \rc+\Lambda_{st} \big(\lc \si(y)^\sharp-\si(\yti)^\sharp \rc \cdot \der x+\si(\yti)^\sharp \cdot \der (x-\xti) \\
+\der \lc y^x \si'(y)-\yti^x \si'(\yti) \rc \cdot \mathbf{x}^2_{st}+\der(\yti^x \si'(\yti)) \cdot \lc \mathbf{x}^2-\tilde{\mathbf{x}}^2 \rc \big),
\end{multline*}
where we have used 
$$(\der y)_{st}= \lc (\operatorname{id}-\Lambda \delta) (\sigma(y) \cdot \der x + (\sigma(y))^x \cdot   \xd) \rc_{st}, $$ some  standard computations yield
\begin{multline*}
\cn[y-\yti;\cq_\ka^{x,\xti}([0,T_0])]\\
 \leq c_{x,\xti,y,\yti} \lcl T_0^\ka \cn[y-\yti;\cq_\ka^{x,\xti}([0,T_0])]+\norm{x-\xti}_\ga+\norm{\mathbf{x}^2-\tilde{\mathbf{x}}^2}_{2\ga}+\lln a-\tilde{a}\rrn \rcl
\end{multline*}
with
$$c_{x,\xti,y,\yti}=c  \lcl 1+\norm{x}_\ga+\norm{\mathbf{x}^2}_{2\ga}+\norm{\xti}_\ga+\norm{\tilde{\mathbf{x}}^2}_{2 \ga}+\cn[y;\cq_\ka^x([0,T])]^2+\cn[\yti;\cq_\ka^{\tilde{x}}([0,T])]^2 \rcl,$$
for some constant $c>0$. Now remember that $\cn[y;\cq_\ka^x([0,T])] \leq P_T(\norm{x}_\ga,\norm{\mathbf{x}^2}_{2\ga})$, as well as $\cn[\yti;\cq_\ka^{\xti}([0,T])] \leq P_T(\norm{\xti}_\ga,\norm{\tilde{\mathbf{x}}^2}_{2\ga})$, for a certain polynomial function $P_T$, so that
$c_{x,\xti,y,\yti} \leq c_{x,\xti}$, where $c_{x,\xti}>0$ stands for a polynomial expression of $\norm{x}_\ga,\norm{\mathbf{x}^2}_{2\ga}$ and $\norm{\xti}_\ga,\norm{\tilde{\mathbf{x}}^2}_{2\ga}$. Set $T_0=(2c_{x,\xti})^{-1/\ka}$ and in this way
$$\cn[y-\yti;\cq_\ka^{x,\xti}([0,T_0])] \leq 2c_{x,\xti} \lcl \norm{x-\xti}_\ga+\norm{\mathbf{x}^2-\tilde{\mathbf{x}}^2}_{2\ga}+\lln a-\tilde{a}\rrn \rcl.$$

\smallskip

\noindent
\textit{Extending the inequality}. With the same arguments as in the above step, we get, for any $k\geq 1$,
\bean
\lefteqn{\cn[y-\yti;\cq_\ka^{x,\xti}([kT_0,(k+1)T_0])]}\\
& \leq & 2c_{x,\xti} \lcl \norm{x-\xti}_\ga+\norm{\mathbf{x}^2-\tilde{\mathbf{x}}^2}_{2\ga}+\lln y_{kT_0}-\yti_{kT_0}\rrn \rcl\\
& \leq & 2c_{x,\xti} \lcl \norm{x-\xti}_\ga+\norm{\mathbf{x}^2-\tilde{\mathbf{x}}^2}_{2\ga}+\lln a-\tilde{a}\rrn +T_0^\ka \sum_{l=0}^{k-1} \cn[y-\yti;\cq_\ka^{x,\xti}([lT_0,(l+1)T_0])]\rcl
\eean
and as a result 
$$\cn[y-\yti;\cq_\ka^{x,\xti}([kT_0,(k+1)T_0])] \leq 2c_{x,\xti}\cdot e^k \lcl \norm{x-\xti}_\ga+\norm{\mathbf{x}^2-\tilde{\mathbf{x}}^2}_{2\ga}+\lln a-\tilde{a}\rrn\rcl$$
using the discrete version of Gronwall's Lemma.

Inequality (\ref{lipsch-cont}) is then a direct consequence of
$$\cn[y-\yti;\cac_1^\ga([0,T])] \leq \sum_{k=0}^{N_{T_0}-1} \cn[y-\yti;\cq_\ka^{x,\xti}([kT_0,(k+1)T_0])],$$
where $N_{T_0}$ is the smallest integer such that $N_{T_0} \cdot T_0 \geq T$, so that $N_{T_0} \leq 1+T/T_0 \leq 1+T \cdot (2c_{x,\xti})^\ka$.

\end{document}